\documentclass[10pt]{article}
 
\usepackage{amsmath}
\usepackage{amsfonts}
\usepackage{hyperref}
\usepackage{graphicx}
\usepackage{multirow}
\usepackage{xcolor}
\usepackage{enumitem}
\usepackage{bibspacing}

\newtheorem{lem}{Lemma}[section]
\newtheorem{teo}{Theorem}[section]
\newtheorem{coro}{Corollary}[section]
\newtheorem{assump}{Assumption}

\newcommand{\halmos}{\hfill$\Box$}
\newenvironment{pro}{\noindent\textit{Proof:}}{\halmos}

\newcommand{\R}{\mathbb{R}}
\newcommand{\N}{\mathbb{N}}
\newcommand{\half}{\frac{1}{2}}
\newcommand{\kbar}{\underline{k}}
\newcommand{\sigmin}{\sigma_{\min}}
\newcommand{\sigmax}{\sigma_{\max}}
\newcommand{\taumin}{\tau_{\min}}
\newcommand{\taumax}{\tau_{\max}}
\newcommand{\spg}{\mathrm{spg}}
\newcommand{\acc}{\mathrm{accel}}
\newcommand{\new}{\mathrm{new}}
\newcommand{\trial}{{\mathrm{trial}}}
\newcommand{\extra}{{\mathrm{extra}}}
\newcommand{\sacc}{s_{\acc}}

\newcommand{\yacc}{y_{\acc}}
\newcommand{\sextra}{s_{\extra}}
\newcommand{\xextra}{x_{\extra}}
\newcommand{\yextra}{y_{\extra}}

\begin{document}

\title{Secant acceleration of sequential residual methods for solving\\
  large-scale nonlinear systems of equations\thanks{This work was
    supported by FAPESP (grants 2013/07375-0, 2016/01860-1, and
    2018/24293-0) and CNPq (grants 302538/2019-4 and 302682/2019-8).}}

\author{
  E. G. Birgin\thanks{Department of Computer Science, Institute of
    Mathematics and Statistics, University of S\~ao Paulo, Rua do
    Mat\~ao, 1010, Cidade Universit\'aria, 05508-090, S\~ao Paulo, SP,
    Brazil. e-mail: egbirgin@ime.usp.br}
  \and
  J. M. Mart\'{\i}nez\thanks{Department of Applied Mathematics,
    Institute of Mathematics, Statistics, and Scientific Computing
    (IMECC), State University of Campinas, 13083-859 Campinas SP,
    Brazil. e-mail: martinez@ime.unicamp.br}}

\date{December 24, 2020\footnote{Typo corrected on January 8,
    2020. Revision made on July 26, 2021.}}

\maketitle

\begin{abstract}
Sequential Residual Methods try to solve nonlinear systems of
equations $F(x)=0$ by iteratively updating the current approximate
solution along a residual-related direction. Therefore, memory
requirements are minimal and, consequently, these methods are
attractive for solving large-scale nonlinear systems. However, the
convergence of these algorithms may be slow in critical cases;
therefore, acceleration procedures are welcome. In this paper, we
suggest to employ a variation of the Sequential Secant Method in order
to accelerate Sequential Residual Methods. The performance of the
resulting algorithm is illustrated by applying it to the solution of
very large problems coming from the discretization of partial
differential equations.\\


\end{abstract}

\section{Introduction} \label{intro}

In the process of solving many real-life problems, it is necessary to
handle large-scale nonlinear systems of equations.  The most obvious
choice for solving these systems is Newton's method, which requires to
solve a possibly large and sparse linear system of equations at each
iteration. Although being very effective in many cases, Newton's
method cannot be employed for solving very large problems when the
Jacobian is unavailable or when it has an unfriendly structure that
makes its factorization unaffordable. On the other hand, Inexact
Newton methods, that solve the Newtonian linear system approximately
at each iteration, are usually
effective~\cite{des,ew1,ew2}. Inexact-Newton methods based on linear
iterative solvers as GMRES may need many matrix-vector products per
iteration. Usually, matrix-vector products of the form $J(x) v $ are
replaced with incremental quotients $[ F(x+hv)-F(x) ] / h$, a
procedure that does not deteriorate the overall performance of
GMRES~\cite{bww,pernice}. However, when GMRES requires many
matrix-vector products for providing a suitable approximate solution
to the Newtonian linear system, the number of residual evaluations per
inexact-Newton iteration may be big. Additional residual evaluations
may also be necessary to decide acceptance of trial points at every
iteration.

This state of facts led to the introduction of algorithms in which the
number of residual evaluations used to compute trial points at each
iteration is minimal, as well as the memory used to store directions
and the computer effort of linear algebra calculations.
DF-SANE~\cite{lmr} was introduced for solving large problems and was
used for solving equilibrium models for the determination of
industrial prices \cite{miller}, multifractal analysis of spot prices
\cite{rypdal}, elastoplastic contact problems \cite{frerot,frerot2},
and PDE equations in reservoir simulations \cite{nelio}, among others.
Improved versions of DF-SANE were given in
\cite{lacruz,meli}. However, the pure form of DF-SANE may be
ineffective for some large problems in which it is necessary to
perform many backtrackings per iteration in order to obtain sufficient
descent. Therefore, on the one hand, it is necessary to investigate
alternative choices of trial steps and, on the other hand,
acceleration procedures are welcome.

Acceleration devices for iterative algorithms are frequent in the
Numerical Analysis
literature~\cite{anderson,brezinski,brezinskiredivo,hoow,niwalker,wwy}. They
incorporate useful information from previous iterations instead of
expending evaluations at the current one. In particular, Anderson's
acceleration introduced in \cite{anderson} is known to produce very
good results when associated with fixed-point
iterations~\cite{bian,chen,fangsaad,hoow,niwalker}, specifically those
originated in Self-Consistent Field (SCF) approaches for electronic
structure calculations \cite{lebris,pulay}. Anderson's acceleration is
closely related to some quasi-Newton methods
\cite{bk,degroote,haelte,fangsaad,lindner} and multipoint secant
algorithms~\cite{barnes,graggstewart,jankowska,martinezbit,or,wolfe}.
The recent survey \cite{brs} sheds a lot of light on the properties of
Anderson's acceleration, generalizations, and relations with other
procedures for accelerating the convergence of sequences.

This work introduces a generalized and accelerated version of
DF-SANE. The generalization consists of allowing non-residual
(although residual-related) directions. The acceleration is based on
the multipoint secant idea, taking advantage of the residual-like
direction steps. Global convergence results that extend the theory of
DF-SANE are given.

The paper is organized as follows. Section~\ref{algorithm} introduces
the accelerated sequential residual methods. Global convergence is
established in Section~\ref{global}. Section~\ref{acceleration}
describes the acceleration process in detail. Implementation features
and numerical experiments are given in Sections~\ref{implementation}
and~\ref{experiments}, respectively. The last section presents the
conclusions.\\

\noindent
\textbf{Notation.} The symbol $\| \cdot \|$ denotes the Euclidean
norm. $\N = \{0,1,2,\dots\}$ denotes the set of natural numbers.
$J(x)$ denotes the Jacobian matrix of $F:\R^n \to \R^n$ computed at
$x$.  For all $x\in\R^n$, we denote $g(x)= J(x)^T F(x)=\nabla \half
\|F(x)\|_2^2$.  If $\{z_k\}_{k \in \N}$ is a sequence and $K = \{k_1,
k_2, k_3, \dots\}$ is an infinite sequence of natural numbers such
that $k_i < k_j$ if $i < j$, we denote
\[
\lim_{k \in K} z_k = \lim_{j \to \infty} z_{k_j}.
\]


\section{Accelerated sequential residual methods} \label{algorithm}

Given $F:\R^n \to \R^n$, consider the problem of finding $x \in \R^n$
such that
\begin{equation} \label{theproblem}
F(x) = 0.
\end{equation}
A radical iterative approach for solving~(\ref{theproblem}) is to
employ only residuals as search directions. Given $\sigma > 0$,
problem~(\ref{theproblem}) is clearly equivalent to $ x = x - \sigma
F(x)$. This trivial observation motivates the introduction of a
fixed-point method given by $x^{k+1} = x^k - \sigma_k F(x^k)$, where
$\sigma_k$ is defined at every iteration. Methods based on this
approach will be called Sequential Residual Methods (SRM) in the
present paper.

Popular SRM were inspired by the Barzilai-Borwein or spectral choice
for the minimization of
functions~\cite{barzilaiborwein,raydan1,raydan2}. Defining
\begin{equation} \label{sy}
s^k = x^{k+1}-x^k \mbox{  and  } y^k = F(x^{k+1}) - F(x^k),
\end{equation}
algorithms SANE~\cite{lacruzraydan} and DF-SANE~\cite{lmr} compute
\begin{equation} \label{alfa}
\sigma_{k+1} = \|s^k\|^2/(y^k)^T s^k,
\end{equation}
safeguarded in such a way that $|\sigma_{k+1}|$ is bounded and bounded
away from zero.  This formula had been used in the context of
self-scaling variable metric methods for minimization \cite{oren} as
it provides a scale invariant diagonal first approximation of the
Hessian. The choice of $\sigma_{k+1}$ may be justified with the same
arguments that Raydan~\cite{raydan1} employed for the choice of the
Barzilai-Borwein or spectral step in minimization problems. After the
computation of $x^{k+1}$, we consider the (generally unsolvable)
problem of satisfying the secant equation \cite{ds} $B_{k+1} s^k =
y^k$ subject to $B_{k+1} = c I$. This leads to the minimization of
$\|c I s^k - y^k\|^2$, whose solution, if $s^k \neq 0$, is $ c =
(y^k)^T s^k/\|s^k\|^2$.  Therefore, a ``natural'' residual-based
iteration for solving problem~(\ref{theproblem}) could be given by $
x^{k+1} = x^k - \sigma_k F(x^k)$, with $\sigma_0$ arbitrary and
$\sigma_{k+1}$ defined by a safeguarded version of (\ref{alfa}) for
all $k \geq 0$.

However, unlike the case of unconstrained minimization, in which
$F(x^k)$ is a gradient, the direction $d^k = - \sigma_k F(x^k)$ may
not be a descent direction for the natural merit function $f(x)$
defined by
\[
f(x) = \half \|F(x)\|^2 \mbox{ for all } x \in \R^n.
\]
In SANE~\cite{lacruzraydan}, a test is performed in order to verify
whether $F(x^k)$ is a descent direction. If this is the case, since
$\nabla f(x) = J(x)^T F(x)$, we should have
\[
F(x^k)^T J(x^k) F(x^k) < 0.
\]
In order to avoid the employment of derivatives, SANE employs the
approximation
\[
J(x^k) F(x^k) \approx \frac{F(x^k + h F(x^k)) - F(x^k)}{h},
\]
for a small $h > 0$. In this way, the descent test is equivalent to
\[
F(x^k)^T F(x^k + h F(x^k)) < \|F(x^k)\|^2,
\]
which requires an auxiliary functional evaluation per iteration. The
necessity of an auxiliary residual evaluation per iteration in SANE
motivated the introduction of DF-SANE \cite{lmr}. Roughly speaking, in
DF-SANE, one gets descent by starting with the trial point $x^k -
\sigma_k F(x^k)$ and proceeding to a double backtracking scheme along
positive and negative directions, aiming that $\|F(x^{k+1})\|$ be
sufficiently smaller than the maximum value of the residual norm in
$M$ consecutive past iterations, where $M$ is given.

The description of the SRM algorithm provided in this section aims to
emphasize the aspects that influence theoretical convergence
properties. For this reason, acceleration steps appear only as a small
detail in the description of the algorithm, although, in practice,
they are essential for the algorithm robustness and efficiency. The
description of the algorithm follows.\\

\noindent
\textbf{Algorithm~\ref{algorithm}.1.} Let $\gamma \in (0,1)$, $0 <
\sigmin < \sigmax < \infty$, $0 < \taumin < \taumax < 1$, a positive
integer~$M$, a sequence $\{\eta_k\}$ such that $\eta_k > 0$ for all $k
\in \N$ and
\begin{equation} \label{sumaetak}
\lim_{k \to \infty} \eta_k = 0,
\end{equation}
and $x_0 \in \R^n$ be given. Set $k \leftarrow 0$.

\begin{description}
\item[Step 1.] If $F(x^k) = 0$, then terminate the execution of the
  algorithm.

\item[Step 2.] Choose $\sigma_k$ such that $|\sigma_k| \in
  [\sigma_{\min}, \sigma_{\max}]$ and $v^k \in \R^n$ such that $\|v^k\|
  = \|F(x^k)\|$. Compute
  \begin{equation} \label{barfk}
    \bar f_k = \max\{f(x^k), \dots, f(x^{\max\{0,k-M+1\}}) \}.
  \end{equation}
  \begin{description}
    \item[Step 2.1.] Set $\alpha_+ \leftarrow 1$ and $\alpha_-
      \leftarrow 1$.
    \item[Step 2.2.] Set $d \leftarrow -\sigma_k v^k$ and $\alpha
      \leftarrow \alpha_+$. Consider
      \begin{equation} \label{algarmijo}
        f(x^k + \alpha d) \leq \bar f_k + \eta_k - \gamma \alpha^2 f(x^k).
      \end{equation}
      If (\ref{algarmijo}) holds, then define $d^k = d$ and $\alpha_k
      = \alpha$ and go to Step~3.
    \item[Step 2.3.] Set $d \leftarrow \sigma_k v^k$ and $\alpha
      \leftarrow \alpha_-$. If (\ref{algarmijo}) holds, then define
      $d^k = d$ and $\alpha_k = \alpha$ and go to Step~3.
    \item[Step 2.4.] Choose $\alpha_+^{\new} \in [\taumin \alpha_+,
      \taumax \alpha_+]$ and $\alpha_-^{\new} \in [\taumin \alpha_-,
      \taumax \alpha_{-}]$, set $\alpha_+ \leftarrow \alpha_+^{\new}$,
      $\alpha_- \leftarrow \alpha_-^{\new}$, and go to Step 2.2.
  \end{description}
  
\item[Step 3.] Compute $x^{k+1}$ such that $f(x^{k+1}) \leq f(x^k +
  \alpha_k d^k)$, set $k \leftarrow k+1$, and go to Step~1.
\end{description}

As in DF-SANE, the sufficient decrease in~(\ref{algarmijo})
corresponds to a nonmonotone strategy that combines the ones
introduced in~\cite{gll} and~\cite{lifuku}. The main differences of
Algorithm~\ref{algorithm}.1 with respect to DF-SANE are the presence
of accelerations at Step~3 and the choice of the non-accelerated step
in a residual-related way, but not necessarily in the residual
direction. The DF-SANE method presented in \cite{lmr} is the
particular case of Algorithm~\ref{algorithm}.1 in which $v^k = F(x^k)$
and $x^{k+1} = x^k + \alpha_k d^k$. The generalization of using a
direction~$v^k$ that satisfies $\|v^k\|=\|F(x^k)\|$ but that is not
necessarily the residue opens up the possibility of using other
alternatives. Moreover, as it will be shown in the next section,
taking random residual-related directions an infinite number of times
makes the method to converge to points at which the gradient of the
sum of squares vanishes. (This property is absent in the original
version of DF-SANE, unless additional conditions on the problem are
required.)  The condition of choosing $x^{k+1}$ that satisfies
$f(x^{k+1}) \leq f(x^k + \alpha_k d^k)$ evidently allows one to choose
$x^{k+1} = x^k + \alpha_k d^k$ as in the case of DF-SANE, but it also
opens the possibility of replacing $x^k + \alpha_k d^k$ by something
even better, precisely opening the possibility for acceleration.

\section{Global convergence} \label{global}

In this section we prove global convergence properties of
Algorithm~\ref{algorithm}.1. Our main purpose is to find solutions of
$F(x)=0$ or, at least, points at which the residual norm $\|F(x)\|$ is
as small as desired, given an arbitrary tolerance. However, this
purpose could be excessively ambitious because, in the worst case,
solutions of the system, or even approximate solutions, may not exist.
For this reason we analyze the situations in which convergence to a
(stationary) point, at which the gradient of the sum of squares
vanishes, occur.  The option of taking directions that are not
residuals but are residual-related in the sense that their norms
coincide with those of the residuals is crucial for this
purpose. Roughly speaking, we will prove that, taking random
residual-related directions an infinite number of times (but not at
all iterations) convergence to stationary points necessarily takes
place.

In Lemma~\ref{lem31} we prove that, given an arbitrary tolerance
$\varepsilon > 0$, either Algorithm~\ref{algorithm}.1 finds an
approximate solution such that $\|F(x^k)\| \leq \varepsilon$ in a
finite number of iterations or produces an infinite sequence of steps
$\{\alpha_k\}$ that tends to zero. For this purpose we will only use
continuity of~$F$.  The lemma begins with a simple proof that the
iteration is well defined.

\begin{lem} \label{lem31}
Assume that $F$ is continuous, $x^k \in \R^n$ is an arbitrary iterate
of Algorithm~\ref{algorithm}.1, and $\|F(x^k)\|\neq 0$. Then,
$\alpha_k$ and $d^k$ satisfying~(\ref{algarmijo}) and $x^{k+1}$
satisfying $f(x^{k+1}) \leq f(x^k + \alpha_k d^k)$ are well
defined. Moreover, if the algorithm generates an infinite sequence
$\{x^k\}$, then there exists an infinite subset of indices $K_1
\subset \N$ such that
\begin{equation} \label{limfk1}
\lim_{k \in K_1} \|F(x^k)\| = 0
\end{equation}
or there exists an infinite subset of indices $K_2$ such that
\begin{equation} \label{lasufi} 
\lim_{k \in K_2} \alpha_k  = 0.
\end{equation}
\end{lem}

\begin{pro}
The fact that $\alpha_k$, $d^k$, and $x^{k+1}$ are always well defined
follows from the continuity of $F$ using that $\eta_k > 0$ and that
$\alpha_k$ can be as small as needed. Assume that the algorithm
generates an infinite sequence $\{x^k\}$. Suppose that~(\ref{limfk1})
does not hold.  Then, there exists $c_1 > 0$ such that
\begin{equation} \label{mayquec}
\|F(x^k)\| > c_1 \mbox{ for all } k \mbox{ large enough}.
\end{equation}  
If (\ref{lasufi}) does not hold either, then there exists 
 $c_2 > 0$ such that
\[
\alpha_k > c_2  \mbox{ for all } k \mbox{ large enough}.
\]
Then, by (\ref{mayquec}),
\begin{equation} \label{nonox}
\alpha_k^2 f(x^k) > (c_1 c_2)^2/2 \mbox{ for all } k \mbox{ large enough}.
\end{equation} 
Then, since $x^{k+1}$ is well defined, for all $k$ large enough,
\begin{equation} \label{sincefa}
f(x^{k+1}) \leq \bar f_k + \eta_k - \gamma \alpha_k^2 f(x^k) \leq \bar f_k +
\eta_k - \gamma (c_1 c_2)^2/2,
\end{equation}
where 
\begin{equation} \label{efek}
\bar f_k = \max\{f(x^k), \dots, f(x^{\max\{0,k-M+1\}}) \}.
\end{equation}        
Since $\eta_k \to 0$, then there exists $k_1 \in \N$ such that for
all $k \geq k_1$,
\begin{equation} \label{sincefa2}
f(x^{k+1}) \leq \bar f_k - \gamma (c_1 c_2)^2/4
\end{equation} 
and, by induction,
 \begin{equation} \label{barmen}
f(x^{k+j}) \leq \bar f_k - \gamma (c_1 c_2)^2/4
\end{equation}  
 for all $j = 1, \dots, M$. Therefore, 
\[
\bar f_{k+M} \leq \bar f_k - \gamma (c_1 c_2)^2/4.
\]
Since this inequality holds for all $k \geq k_1$ we conclude that
$\bar f_{k+jM}$ is negative for $j$ large enough, which contradicts
the fact that $f(x) \geq 0$ for all $x \in \R^n$.
\end{pro}

\begin{lem} \label{teo31}
Assume that $F(x)$ admits continuous derivatives for all $x \in \R^n$
and $\{x^k\}$ is generated by Algorithm~\ref{algorithm}.1. Assume that
$\{\|F(x^k)\|\}$ is bounded away from zero and that $K_2$ is an
infinite set of indices such that
 \begin{equation} \label{alfaacero}
\lim_{k \in K_2} \alpha_k = 0.
\end{equation}        
 Let $x_* \in \R^n$ be a limit point of $\{x^k\}_{k \in K_2}$. Then,
 the set of limit points of $\{v^k\}_{k \in K_2}$ is
 nonempty; and, if $v$ is a limit point of $\{v^k\}_{k \in K_2}$,
 then we have that
\begin{equation} \label{ortogo}
\langle J(x_*) v, F(x_*) \rangle = \langle v, J(x_*)^T F(x_*) \rangle = 0.
\end{equation}
\end{lem}

\begin{pro}
Assume that $K \subset K_2$ is such that $\lim_{k \in K} x^k =
x_*$. Then, by continuity, $\lim_{k \in K} F(x^k) =
F(x_*)$. Therefore, the sequence $\{F(x^k)\}_{k \in K}$ is
bounded. Since $\|v^k\| = \|F(x^k)\|$ for all~$k$, we have that the
sequence $\{v^k\}_{k \in K}$ is bounded too; therefore, it admits a
limit point $v \in \R^n$, as we wanted to prove.

Let $K_3$ be an infinite subset of $K$ and $v \in \R^n$ be such that
$\lim_{k \in K_3} v^k = v$. Since the first trial value for $\alpha_k$
at each iteration is $1$, (\ref{alfaacero}) implies that there exist
$K_4$, an infinite subset of $K_3$, $\alpha_{k,+} > 0$, $\alpha_{k,-}
> 0$, $\sigma_k \in [\sigmin, \sigmax]$, and $d^k = -\sigma_k v^k$
such that
\[
\lim_{k \in K_4} \sigma_k = \sigma \in   [\sigmin, \sigmax],
\]
\[
\lim_{k \in K_4} \alpha_{k,+} = \lim_{k \in K_4} \alpha_{k,-} = 0, 
\]
and, for all $k \in K_4$, 
\[
f(x^k + \alpha_{k,+} d^k) >  \bar f_k + \eta_k - \gamma \alpha_{k,+}^2 f(x^k), 
\]
and 
\[
f(x^k - \alpha_{k,-} d^k) > \bar f_k + \eta_k - \gamma \alpha_{k,-}^2 f(x^k).
\]
Therefore, by the definition~(\ref{barfk}) of $\bar f_k$, for all $k
\in K_4$,
\[
f(x^k + \alpha_{k,+} d^k) >  f(x^k) + \eta_k - \gamma \alpha_{k,+}^2 f(x^k)
\]
and 
\[
f(x^k - \alpha_{k,-} d^k) > f(x^k) + \eta_k - \gamma \alpha_{k,-}^2 f(x^k).
\]
So, since $\eta_k > 0$, 
\[
\frac{f(x^k + \alpha_{k,+} d^k) - f(x^k)}{\alpha_{k,+}} > - \gamma
\alpha_{k,+} f(x^k)
\]
and 
\[
\frac{ f(x^k - \alpha_{k,-} d^k) - f(x^k)}{\alpha_{k,-}} > - \gamma
\alpha_{k,-} f(x^k)
\]
for all $k \in K_4$. Thus, by the Mean Value Theorem, there exist
$\xi_{k,+} \in [0, \alpha_{k,+}]$ and $\xi_{k,-} \in [0,
  \alpha_{k,-}]$ such that
\begin{equation} \label{notarmijo7}
\langle \nabla f(x^k + \xi_{k,+} d^k), d^k \rangle > - \gamma
\alpha_{k,+} f(x^k)
\end{equation}
and 
\begin{equation} \label{notarmijo8}
-\langle \nabla f(x^k - \xi_{k,-} d^k), d^k \rangle > - \gamma
\alpha_{k,-} f(x^k)
\end{equation}
for all $k \in K_4$. Taking limits for $k \in K_4$ in both sides of
(\ref{notarmijo7}) and (\ref{notarmijo8}) we get
\begin{equation} \label{notarmijo9}       
\langle \nabla f(x_*), \sigma v \rangle = 0.
\end{equation}
Therefore, (\ref{ortogo}) is proved. 
\end{pro}\\

Theorem~\ref{teo32} will state a sufficient condition for the
annihilation of the gradient of $f(x)$ at a limit point of the
sequence. For that purpose, we will assume that the absolute value of
the cosine of the angle determined by $v^k$ and $J(x^k)^T F(x^k)$ is
bigger than a tolerance $\omega > 0$ infinitely many times. Note that
we do not require this condition to hold for every $k \in \N$. The
distinction between these two alternatives is not negligible. For
example, if $v^k$ is chosen infinitely many times as a random vector
the required angle condition holds with probability~1. Conversely, if
we require the fulfillment of the angle condition for all $k \in \N$
and we choose random directions, the probability of fulfillment would
be zero.

\begin{teo} \label{teo32}
Assume that $F(x)$ admits continuous derivatives for all $x \in \R^n$
and $\{x^k\}$ is generated by Algorithm~\ref{algorithm}.1. Suppose
that $\sum_{k=0}^\infty \eta_k = \eta < \infty$ and the level set $\{x
\in \R^n \;|\; f(x) \leq f(x^0) + \eta\}$ is bounded. Assume that at
least one of the following possibilities hold:
\begin{enumerate} 

\item Agorithm~\ref{algorithm}.1 stops at Step~1 for some $k \in \N$.

\item There exists a sequence of indices $K_1 \subset \N$ such that
  $\lim_{k \in K_1} \|F(x^k)\|=0$.

\item There exist $\omega > 0$ and a sequence of indices $K_2$ such
  that
  \[
  \lim_{k \in K_2} \alpha_k = 0
  \]
  and  
  \begin{equation} \label{cosangle}
    |\langle v^k, J(x^k)^T F(x^k) \rangle| \geq \omega \|v^k\|
    \|J(x^k)^T F(x^k)\|
  \end{equation}
  for all $k \in K_2$. 
\end{enumerate}
Then, for any given $\varepsilon > 0$, there exists $k \in \N$ such
that $\|J(x^k)^T F(x^k)\| \leq \varepsilon$.  Moreover, if the
algorithm does not stop at Step~1 with $\|F(x^k)\|=0$, there exists a
limit point $x_*$ of the sequence generated by the algorithm such that
$\|J(x_*)^T F(x_*)\| = 0$.
\end{teo}

\begin{pro}
If the sequence generated by the algorithm stops at some $k$ such that
$\|F(x^k)\|=0$ the thesis obviously holds.

By the boundedness of the level set defined by $f(x^0)+\eta$, the
sequence $\{x^k\}$ is bounded.  Therefore, by continuity,
$\{\|F(x^k)\|\}$ and $\{\|J(x^k)\|\}$ are also bounded and $\lim_{k
  \in K_1} \|F(x^k)\|=0$ implies that $\lim_{k \in K_1} \|J(x^k)^T
F(x^k)\| =0$.  So, the thesis also holds when the second alternative
in the hypothesis takes place.

Therefore, we only need to consider the case in which $\|F(x^k)\|$ is
bounded away from zero and the third alternative in the hypothesis
holds.  By the definition of $v^k$ and the boundedness of $\|F(x^k)\|$
we have that $\{\|v^k\|\}$ is also bounded.  Therefore, there exist an
infinite set $K_3 \subset K_2 $, $x_* \in \R^n$, and $v \in \R^n$ such
that $\lim_{k \in K_3} x^k = x_*$ and $\lim_{k \in K_3} v^k =
v$. Therefore, by Lemma~\ref{teo31},
\[
\langle v, J(x_*)^T F(x_*) \rangle = 0.
\]
So, by the convergence of $\{v^k\}$ and $\{x^k\}$ for $k \in K_3$ and
the continuity of $F$ and $J$,
\[
\lim_{k \in K_3} \langle v^k, J(x^k)^T F(x^k) \rangle = 0.
\]
Therefore, by (\ref{cosangle}), since $K_3 \subset K_2$,
\begin{equation} \label{fk1}
\lim_{k \in K_3}  \|v^k\| \|J(x^k)^T F(x^k)\|   = 0.
\end{equation} 
Since, in this case, $\|F(x^k)\|$ is bounded away from zero, we have
that $\|v^k\|$ is bounded away from zero too. Therefore, by
(\ref{fk1}),
\begin{equation} \label{fk11}
\lim_{k \in K_3} \|J(x^k)^T F(x^k)\| = 0.
\end{equation}    
Therefore, the thesis is proved.
\end{pro} \\

\noindent
\textbf{Remark.} The hypotheses of Theorem~\ref{teo32} are not
plausible for the classical DF-SANE algorithm because, in general, it
is impossible to guarantee the angle condition~(\ref{cosangle}) if we
only use residual directions. In that case, convergence to points such
that $F(x^*)$ is orthogonal to $J(x^*)^T F(x^*)$ but $\|J(x^*)^T
F(x^*)\| \neq 0$ is possible. The generalization introduced in this
paper for the choice of $v^k$ provides a remedy for that drawback. In
order to fulfill the requirements of Theorem~\ref{teo32},
Algorithm~\ref{algorithm}.1 can be implemented in the following way:
\begin{itemize}
\item At the beginning of the algorithm, we define
  $\alpha_{\mathrm{small}} > 0$.
\item At every iteration of the algorithm, if $\alpha_k <
  \alpha_{\mathrm{small}}$, then we discard the direction, we
  define $v^k/\|v^k\|$ as a random point, with uniform distribution,
  in the unitary sphere, and we complete the computation of~$x^{k+1}$
  executing Steps~2.1--2.4 and Step~3 of the algorithm. If the new
  computed $\alpha_k$ is smaller than $\alpha_{\mathrm{small}}$, then
  we update $\alpha_{\mathrm{small}} \leftarrow
  \alpha_{\mathrm{small}}/2$.
\end{itemize}
Note that, if $\alpha_{\mathrm{small}}$ is updated only a finite
number of times, we have that $\alpha_k$ turns out to be bounded away
from zero and, by Lemma~\ref{lem31}, $\|F(x^k)\|$ tends to zero.
Therefore, we only need to consider the case in which
$\alpha_{\mathrm{small}}$ is updated infinitely many times. In this
case we have that $\alpha_k$ tends to zero along the subsequence that
corresponds to $\alpha_{\mathrm{small}}$ updates and, consequently, the
random choice of $v^k$ is made for a whole subsequence along which
$\alpha_k$ tends to zero. For each random choice the probability
that~(\ref{cosangle}) does not hold is strictly smaller than~$1$.
Therefore, the probability of existence of a subsequence such that all
its members satisfy~(\ref{cosangle}) is~$1$.

\begin{coro} \label{coro30}
Suppose that the assumptions of Theorem~\ref{teo32} hold. Assume,
moreover, that there exists $\gamma > 0$ such that for all $k \in \N$,
\begin{equation} \label{gama}
\|J(x^k)^T F(x^k)\| \geq \gamma \|F(x^k)\|.
\end{equation}
Then, given $\varepsilon > 0$, there exists an iterate $k$ such that
$\|F(x^k)\| \leq \varepsilon$.  Moreover, there exists a limit point
$x_*$ of the sequence generated by the algorithm such that $\|F(x_*)\|
= 0$.
\end{coro}   

\begin{pro}
The thesis of this corollary follows directly from (\ref{gama}) and
Theorem~\ref{teo32}.
\end{pro}

\begin{coro} \label{coro40}
Suppose that the assumptions of Theorem~\ref{teo32} hold and $J(x)$ is
nonsingular with $\|J(x)^{-1}\|$ uniformly bounded for all $x \in
\R^n$. Then, for all $\ell = 1, 2, \dots$, there exists an iterate
$k(\ell)$ such that $\|F(x^{k(\ell)})\| \leq 1/\ell$.  Moreover, at
every limit point $x_*$ of the sequence $\{x^{k(\ell)}\}$ we have that
$\|F(x_*)\| = 0$.
\end{coro}

 
The lemma below shows that, if $\| F(x^k) \| \to 0$ for a subsequence,
then it goes to zero for the whole sequence.  We will need the
following Assumptions~\ref{assumpA} and
\ref{assumpB}. Assumption~\ref{assumpA} guarantees that the difference
between consecutive iterates is not greater than a multiple of the
residual norm. Note that such assumption holds trivially with $c =
\sigma_{\max}$ when the iteration is not accelerated and, so, $x^{k+1}
= x^k + \alpha_k d^k$. Therefore, Assumption~\ref{assumpA} needs to be
stated only with respect to accelerated steps; and in practice can be
granted by simply discarding the accelerated iterates that do not
satisfy it. Assumption~\ref{assumpB} merely states that $F$ is
uniformly continuous at least restricted to the set of
iterates. Assumption~\ref{assumpB} holds if~$F$ is uniformly
continuous onto a level set of~$f$ that contains the whole sequence
generated by the algorithm.  A sufficient condition for this fact is
the fulfillment of a Lipschitz condition by the function~$F$.

\begin{assump} \label{assumpA}
There exists $c> 0$ such that, for all $k \in \N$, whenever $x^{k+1}
\neq x^k + \alpha_k d^k$, we have that
\[
\|x^{k+1} -  x^k \| \leq c \|F(x^k)\|.
\]
\end{assump}

\begin{assump} \label{assumpB}
For all $\varepsilon > 0$, there exists $\delta$ such that, whenever
$\|x^{k+1}-x^k\| \leq \delta$ one has that $\|F(x^{k+1})-F(x^k)\| \leq
\varepsilon$.
\end{assump}       

Lemma~\ref{lem32} below is stronger than Theorem~2 of \cite{lmr}
because here we do not assume the existence of a limit point.

\begin{lem} \label{lem32}
Suppose that $F$ is continuous, $\sum_{k = 0}^\infty \eta_k < \infty$,
Algorithm~\ref{algorithm}.1 generates the infinite sequence $\{x^k\}$,
Assumptions~\ref{assumpA} and \ref{assumpB} hold, and there exists an
infinite subsequence defined by the indices in $K_1$ such that
\begin{equation} \label{notenianumero}
\lim_{k \in K_1} \|F(x^k)\| = 0.
\end{equation}
Then,
\begin{equation} \label{limefe}
\lim_{k \to \infty} \|F(x^k)\|=0
\end{equation}
and
\begin{equation} \label{limequis}
\lim_{k \to \infty} \|x^{k+1}-x^k\| = 0.
\end{equation}
\end{lem}

\begin{pro}
Assume that (\ref{limefe}) is not true. Then, there exists an infinite
subsequence $\{x^k\}_{k \in K_2}$ such that
\begin{equation} \label{lak2}
f(x^k) > \hat c > 0
\end{equation}
for all $k \in K_2$. 

Let us prove first that
\begin{equation} \label{kauno}
\lim_{k \in K_1} \| F(x^{k+1}) \|=0.
\end{equation}
By the definition of Algorithm~\ref{algorithm}.1, when $x^{k+1} = x^k
+ \alpha_k d^k$, $\| x^{k+1} - x^k \| \leq \sigma_{\max} \| F(x^k)
\|$. This, along with Assumption~\ref{assumpA} and
(\ref{notenianumero}), implies that
\begin{equation} \label{limequisK1}
\lim_{k \in K_1} \|x^{k+1}-x^k\| = 0.
\end{equation}
By Assumption~\ref{assumpB}, (\ref{limequisK1}) implies
\begin{equation} \label{limFK1}
\lim_{k \in K_1} \|F(x^{k+1})-F(x^k)\| = 0,
\end{equation}
and (\ref{kauno}) follows from (\ref{notenianumero}) and (\ref{limFK1}).


By induction, we deduce that for all $j = 1, \dots, M-1$,
\[
\lim_{k \in K_1} \| F(x^{k+j}) \|=0.
\]
Therefore, 
\begin{equation} \label{aceromevoy}
\lim_{k \in K_1 } \max\{f(x^k), \dots, f(x^{k+M-1})\} = 0.
\end{equation}
But 
\[
\begin{array}{rcl}
f(x^{k+M}) &\leq& \max\{f(x^k), \dots, f(x^{k+M-1})\} + \eta_{k+M-1} -
\gamma \alpha_{k+M-1}^2 f(x^{k+M-1})\\[2mm] &\leq& \max\{f(x^k), \dots,
f(x^{k+M-1})\} + \eta_{k+M-1}.
\end{array}
\]
Analogously,
\[
\begin{array}{rcl}
f(x^{k+M+1}) &\leq& \max\{f(x^{k+1}), \dots, f(x^{k+M})\} + \eta_{k+M} -
\gamma \alpha_{k+M}^2 f(x^{k+M})\\[2mm]
&\leq& \max\{f(x^{k+1}), \dots, f(x^{k+M})\} + \eta_{k+M} \\[2mm]
&\leq& \max\{f(x^{k}), \dots, f(x^{k+M-1})\} + \eta_{k+M-1} + \eta_{k+M};
\end{array}
\]
and, inductively, 
\[
f(x^{k+M+j}) \leq \max\{f(x^{k}), \dots, f(x^{k+M-1})\} + \eta_{k+M-1} +
\eta_{k+M} + \dots + \eta_{k+M+j-1}
\]                   
for all $j = 0, 1, \dots, M-1$. Therefore,
\begin{equation} \label{taking}
\max\{f(x^{k+M}), \dots, f(x^{k+2 M-1})\} \leq \max\{f(x^{k}), \dots,
f(x^{k+M-1})\} + \sum_{j = k+M-1}^{k+2M-2} \eta_j.
\end{equation}
By induction on $\ell = 1, 2, \dots$, we obtain that
\begin{equation} \label{iducele}
\max\{f(x^{k+\ell M}), \dots, f(x^{k+(\ell+1) M-1})\} \leq
\max\{f(x^{k}), \dots, f(x^{k+M-1})\} + \sum_{j = k+M-1}^{k+ (\ell+1)
  M-2} \eta_j.
\end{equation}     
Therefore, for all $\ell = 1, 2, \dots$ we have:
\begin{equation} \label{iduce}
\max\{f(x^{k+\ell M}), \dots, f(x^{k+(\ell+1) M-1})\} \leq
\max\{f(x^{k}), \dots, f(x^{k+M-1})\} + \sum_{j = k+M-1}^{\infty}
\eta_j.
\end{equation}   
  
By (\ref{aceromevoy}), the summability of $\eta_j$ and (\ref{iduce}),
there exists $k_1 \in K_1$ such that for all $k \in K_1$ such that $k
\geq k_1$,
\[
\max\{f(x^{k}), \dots, f(x^{k+M-1})\} + \sum_{j = k+M-1}^{\infty} /
\eta_j \leq \hat c / 2.
\]
Therefore, by (\ref{iduce}), for all $k \in K_1$ such that $k \geq
k_1$ and all $\ell = 1, 2, \dots$,
\begin{equation} \label{idu}
\max\{f(x^{k+\ell M}), \dots, f(x^{k+(\ell+1) M-1})\} \leq \hat c / 2. 
\end{equation}   
Clearly, this is incompatible with the existence of a subsequence such
that (\ref{lak2}) holds. Therefore, (\ref{limefe}) is proved. Then, by
Assumption~\ref{assumpA} and (\ref{limefe}), (\ref{limequis}) also
holds.
\end{pro}

\begin{teo} \label{teofi} 
Suppose that the assumptions of Theorem~\ref{teo32} and
Lemma~\ref{lem32} hold, $J(x)$ is nonsingular and $\|J(x)^{-1}\|$ is
uniformly bounded for all $x \in \R^n$.  Then, there exists $x_* \in
\R^n$ such that $F(x_*) = 0$ and $\lim_{k \to \infty} x^k = x_*$.
\end{teo}

\begin{pro}
By Corollary~\ref{coro40}, there exists a subsequence $\{x^k\}_{k \in
  K_1}$ that converges to a point $x_*$ such that $F(x_*)=0$. Then, by
Lemma~\ref{lem32},
\begin{equation} \label{limff}
\lim_{k \to \infty} F(x^k) = 0
\end{equation}
and
\[
\lim_{k \to \infty} \|x^{k+1}-x^k\|=0.
\]
Since $J(x_*)$ is nonsingular, by the Inverse Function Theorem, there
exists $\delta > 0$ such that $\|F(x)\| > 0$ whenever $0 < \|x - x_*\|
\leq \delta$. Therefore, given $\varepsilon \in (0, \delta]$
  arbitrary, there exists $\hat c > 0$ such that
\[
\|F(x)\| \geq \hat c \mbox{  whenever  } \varepsilon \leq \|x - x_*\| \leq \delta.
\]
By (\ref{limff}), the number of iterates in the set defined by $
\varepsilon \leq \|x - x_*\| \leq \delta$ is finite. On the other
hand, if $k$ is large enough, one has that $\|x^{k+1} - x^k\| \leq
\delta - \varepsilon$. Since there exists a subsequence of $\{x^k\}$
that tends to $x_*$ it turns out that $\|x^k - x_*\| \leq \varepsilon$
if $k$ is large enough. Since $\varepsilon$ was arbitrary this implies
that $\lim_{k \to \infty} x^k = x_*$ as we wanted to prove.
\end{pro}\\

We finish this section with a theorem that states that, under some
assumptions, if in Algorithm~\ref{algorithm}.1 we have $x^{k+1}=x^k +
\alpha_k d^k$ for all~$k$, i.e.\ without accelerations, then $\{x^k\}$
converges superlinearly to a solution. A similar result was proved in
\cite[Thm.~3.1]{fmmr} with respect to a spectral gradient methods with
retards for quadratic unconstrained minimization. We will need two
assumptions. In Assumption~\ref{dennismo}, the sequence generated by
Algorithm~\ref{algorithm}.1 is assumed to be generated by the spectral
residual choice with the first step accepted at each iteration and
without acceleration.  In Assumption~\ref{lipscha}, the Jacobian
$J(x)$ is assumed to be Lipschitz-continuous. The first part of
Assumption~\ref{dennismo} ($v^k = F(x^k)$) limits the consequences of
Theorem~\ref{teo35} to a particular choice at Step~2 of
Algorithm~\ref{algorithm}.1. The fact that $\sigma_k$ is well defined
in the form stated in this assumption and the fact that the descent
condition~(\ref{algarmijo}) holds with $\alpha = 1$ are in fact
assumptions on the behavior of the algorithm, with the same status as
the hypotheses~(\ref{limxk}) and~(\ref{limsk}) of
Theorem~\ref{teo35}. We do not make any claim about the frequency with
which these hypotheses hold.

\begin{assump} \label{dennismo}
Assume that, at Step~2 of Algorithm~\ref{algorithm}.1, we choose
\[
v^k = F(x^k)
\]
and that there exists $k_0 \geq 1$ such that for all $k \geq k_0$, 
\[
(s^{k-1})^T y^{k-1} \neq 0
\]
and 
\[
\sigma_k =   \frac{(s^{k-1})^T s^{k-1}}{(s^{k-1})^T y^{k-1}}, 
\] 
where $s^k$ and $y^k$ are defined by (\ref{sy}).  Assume, moreover,
that, for all $k \geq k_0$, (\ref{algarmijo}) holds with $\alpha = 1$
and, at Step~3, $x^{k+1} = x^k + \alpha_k d^k $.
\end{assump}

By the Mean Value Theorem $y^{k-1} = [\int_0^1 J(x^{k-1} + t s^{k-1})
  dt] s^{k-1}$.  Therefore,
\[
(s^{k-1})^T y^{k-1} = \int_0^1 [(s^{k-1})^T J(x^{k-1} + t s^{k-1})
  s^{k-1}] dt.
\]
Suppose that for all $x \in \R^n$ and $v \neq 0$ we have that $v^T
J(x) v / v^T v $ is bounded and bounded away from zero. This implies
that $(s^{k-1})^T y^{k-1}/(s^{k-1})^T s^{k-1}$ is bounded and bounded
away from zero whenever $s^{k-1} \neq 0$. Analogously, if $- v^T J(x)
v/v^T v $ is bounded and bounded away from zero we have that
$(s^{k-1})^T y^{k-1}/(s^{k-1})^T s^{k-1}$ is bounded and bounded away
from zero. Since the sequence does not terminate at Step~1 we have
that $s^{k-1} \neq 0$ for all $k$. Therefore we see that a sufficient
condition for the existence of $\sigma_{\min}$ and $\sigma_{\max}$
(parameters of Algorithm~\ref{algorithm}.1) such that $|\sigma_k| \in
[\sigma_{\min},\sigma_{\max}]$ holds with $\sigma_k = (s^{k-1})^T
s^{k-1} / ( (s^{k-1})^T y^{k-1} )$ is the boundedness and (positive or
negative) definiteness of the Jacobian $J(x)$.

\begin{assump} \label{lipscha}
There exists $L > 0$ such that for all $x, z $ in an open and convex
set that contains the whole sequence $\{x^k\}$ generated by
Algorithm~\ref{algorithm}.1, the Jacobian is continuous and satisfies
\[
\|J(z)-J(x)\| \leq L \|z - x\|.
\]
\end{assump}

\begin{teo}\label{teo35} 
Assume that the sequence $\{x^k\}$, generated by
Algorithm~\ref{algorithm}.1, does not terminate at Step~1 and that
Assumptions~\ref{dennismo} and~\ref{lipscha} hold.  Suppose that $x_*
\in \R^n$ is such that
\begin{equation} \label{limxk}
\lim_{k \to \infty} x^k = x_*, 
\end{equation}
the Jacobian  $J(x_*)$ is nonsingular 
and $s \in \R^n$ is such that 
\begin{equation} \label{limsk}
\lim_{k\to \infty} \frac{s^k}{\|s^k\|} = s.
\end{equation}
Then, 
\begin{equation} \label{limalfa}
 \lim_{k \to \infty} \frac{(s^k)^T y^k}{(s^k)^T s^k} =  s^T J(x_*) s,  
\end{equation}
\begin{equation} \label{fescero}
F(x_*) = 0, 
\end{equation}
and
\begin{equation} \label{super}
x^k \mbox{ converges $Q$-superlinearly to } x_*.
\end{equation}                          
\end{teo}

\begin{pro}
By Assumption \ref{lipscha}, there exists an open and convex set that
contains the whole sequence $\{x^k\}$ such that for all $x, z$ in that
set,
\begin{equation} \label{lipscho}
F(z) = F(x) + J(x) (z - x) + O(\|z-x\|^2).
\end{equation}
By Assumption~\ref{dennismo}, the sequence $\{x^k\}$ is such that for
all $k \geq k_0$,
\begin{equation} \label{secalfa}
x^{k+1} = x^k - \frac{1}{\kappa_k} F(x^k),
\end{equation}
 where 
\begin{equation} \label{alfak}
\kappa_{k+1} = \frac{(s^k)^T y^k}{(s^k)^T s^k}.
\end{equation}

By (\ref{lipscho}),
\begin{equation} \label{lipscho2}
y^k =  J(x^k) s^k + O(\|s^k\|^2)
\end{equation}
for all $k \geq k_0$. Therefore, by (\ref{lipscho2}),
\begin{equation} 
\kappa_{k+1} =  \frac{(s^k)^T [ J(x^k) s^k + O(\|s^k\|^2) ] }{(s^k)^T s^k}.  
\end{equation}  
Thus, by (\ref{limxk}), (\ref{limsk}), and the continuity of $J(x)$, 
\begin{equation} \label{limalfamas}
\lim_{k \to \infty} \kappa_{k+1} =  s^T J(x_*) s.
\end{equation}
Therefore, (\ref{limalfa}) is proved.

By (\ref{lipscho}) we have that 
\[
F(x^{k+1}) = F(x^k) + J(x^k) s^k + O(\|s^k\|^2).
\]
Then, by (\ref{secalfa}), 
\[
- \kappa_{k+1} s^{k+1} = - \kappa_k s^k   + J(x^k) s^k + O(\|s^k\|^2)
\]   
for all $k \in \N$.  Therefore, for all $k \in \N$,
\[
s^{k+1} = \frac{ \kappa_k s^k - J(x^k) s^k}{\kappa_{k+1}} - \frac{
  O(\|s^k\|^2)}{ \kappa_{k+1}}.
\]    
Therefore, 
\[
s^{k+1} = - \frac{1}{\kappa_{k+1}}[ (J(x^k) - \kappa_k I) s^k +  O(\|s^k\|^2)].
\]   
Then
\[
\frac{s^{k+1}}{\|s^{k+1}\|} = \frac{ - [ (J(x^k) - \kappa_k I) s^k +
    O(\|s^k\|^2)] } { \|- [ (J(x^k) - \kappa_k I) s^k +
    O(\|s^k\|^2)]\|}.
\]
So,
\begin{equation} \label{both}
\frac{s^{k+1}}{\|s^{k+1}\|} = \frac{ - [ (J(x^k) - \kappa_k I)
    s^k/\|s^k\| + O(\|s^k\|^2)/\|s^k\|] } { \|- [ (J(x^k) - \kappa_k
    I) s^k/\|s^k\| + O(\|s^k\|^2)/\|s^k\|]\|}.
\end{equation}   
Define
\begin{equation} \label{alfastar}
\kappa_* =  s^T J(x_*) s.  
\end{equation}   
Assume that 
\begin{equation} \label{ce}
\|(J(x_*) - \kappa_* I)  s \| \neq 0.
\end{equation}
Then, by (\ref{ce}), (\ref{limalfa}), and (\ref{limsk}), taking limits
in both sides of (\ref{both}),
\begin{equation} \label{boz}
s = \frac{ - [ (J(x_*) - \kappa_* I) s ] } { \|- [ (J(x_*) - \kappa_*
    I) s\| }.
\end{equation}                 
Pre-multiplying both sides of (\ref{boz}) by $s^T$, as $s^T s = 1$, we
obtain:
\begin{equation} \label{boz2}
1 = \frac{ - [ (s^T J(x_*) s - \kappa_* ] } { \|- [ (J(x_*) - \kappa_*
    I) s\| }.
\end{equation}
Then, by (\ref{alfastar}), we get a contradiction.  This implies that
the assumption (\ref{ce}) is false. So,
\begin{equation} \label{ceno}
\|(J(x_*) - \kappa_* I)  s \| = 0.
\end{equation}      
Therefore, 
\begin{equation} \label{demor}
\lim_{k \to \infty} \frac{ \kappa_k I s^k - J(x_*) s^k}{\|s^k\|} = 0.
\end{equation}
But (\ref{demor}) is the Dennis-Mor\'e condition \cite{dennismore}
related with the iteration $x^{k+1} = x^k - B_k^{-1} F(x^k)$ with $B_k
= \kappa_k I$. Then, we have that $F(x_*) = 0$ and the convergence is
superlinear.
\end{pro}

\section{Acceleration } \label{acceleration}

The Sequential Secant Method
(SSM)~\cite{barnes,graggstewart,jankowska,martinezbit,wolfe}, called
Secant Method of $n+1$ points in the classical book of Ortega and
Rheinboldt~\cite{or}, is a classical procedure for solving nonlinear
systems of equations. Given $n+1$ consecutive iterates $x^k, x^{k-1},
\dots, x^{k-n}$, the SSM implicitly considers an interpolary affine
function $A: \R^n \to \R^n$ such that $A(x^j) = F(x^j)$ for $j = k,
k-1, \dots, k-n$; and defines $x^{k+1}$ as a solution of $A(x)=0$. If
the points $x^k, x^{k-1}, \dots, x^{k-n}$ are affinely independent
there is only one affine function with the interpolatory property,
which turns out to be defined by:
\[
A(x) = F(x^k) + (y^{k-1}, \dots, y^{k-n})(s^{k-1}, \dots, s^{k-n})^{-1} (x - x^k),
\]
where 
\[
s^j = x^{j+1} - x^j \mbox{ and } y^j = F(x^{j+1}) - F(x^j) \mbox{ for all } j \leq k-1.
\]
Moreover, if $F(x^k), F(x^{k-1}), \dots, F(x^{k-n})$ are also
affinely independent we have that $y^{k-1}, \dots, y^{k-n}$ are
linearly independent and $A(x) = 0$ has a unique solution given by:
\begin{equation} \label{ssmplus}
  x^{k+1} = x^k - (s^{k-1}, \dots, s^{k-n})(y^{k-1}, \dots, y^{k-n})^{-1} F(x^k).
\end{equation}
Practical implementations of SSM try to maintain nonsingularity of the
matrices $(s^{k-1}, \dots, s^{k-n})$ and $(y^{k-1}, \dots$, $y^{k-n})$
employing auxiliary points in the case that linear independence
fails. See \cite{martinezbit} and references therein.
   
If $s^{k-1}, \dots, s^{k-n}$ are linearly independent but $y^{k-1},
\dots, y^{k-n}$ are not, a solution of $A(x) = 0$ may not exist. In
this case we may define $x^{k+1}$ as a solution of
\[
\mbox{Minimize} \| F(x^k) + (y^{k-1}, \dots, y^{k-n})(s^{k-1},
\dots, s^{k-n})^{-1} (x - x^k) \|^2.
\]
A sensible choice for solving this problem is:
\[
x^{k+1} = x^k - (s^{k-1}, \dots, s^{k-n})(y^{k-1}, \dots,
y^{k-n})^{\dagger} F(x^k),
\]     
where $Y^\dagger$ denotes the Moore-Penrose pseudoinverse of a
matrix~$Y$.

This formula suggests a generalization of the SSM. Instead of using
$x^k$ and the $n$ previous iterates to define $x^{k+1}$ we may employ
$p < n$ previous iterates. This leads to the iteration
\begin{equation} \label{ssmp}
  x^{k+1} = x^k - (s^{k-1}, \dots, s^{k-p})(y^{k-1}, \dots, y^{k-p})^{\dagger} F(x^k). 
\end{equation}
Formula (\ref{ssmp}) admits a nice geometrical interpretation. Let
$F_*$ be the vector of minimal norm in the affine subspace spanned by
$F(x^k), \dots, F(x^{k-p})$.  Then,
\begin{equation} \label{efestar}
F_* = F(x^k) + \lambda_{k-1} (F(x^k) - F(x^{k-1})) + \dots +
\lambda_{k-p} (F(x^{k-p+1}) - F(x^{k-p})).
\end{equation}
The idea in (\ref{ssmp}) consists of using the same coefficients
$\lambda_{k-1}, \dots, \lambda_{k-p}$ in order to compute:
\[
x^{k+1} = x^k   + \lambda_{k-1} (x^k - x^{k-1}) + \dots + \lambda_{k-p} (x^{k-p+1} - x^{k-p}). 
\]
In other words, we apply to $x^k, x^{k-1}, \dots, x^{k-p}$ the same
transformation that is applied to $F(x^k), F(x^{k-1}), \dots$,
$F(x^{k-p})$ for obtaining the minimum norm element in the affine
subspace generated by these vectors. The only ambiguity in the
representation above is that the coefficients $\lambda_{k-1}, \dots,
\lambda_{k-p}$ are not unique if $F(x^k), F(x^{k-1}), \dots$,
$F(x^{k-p})$ are not affinely independent. This ambiguity can be
removed taking $\lambda_{k-1}, \dots, \lambda_{k-p}$ that minimize
$\sum_{j=1}^p \lambda_{k-j}^2$ subject to (\ref{efestar}). This
decision leads to the use of the Moore-Penrose pseudoinverse
in~(\ref{ssmp}).

However, although $x^{k+1}$ may represent an improvement with respect
to $x^k$, employing $p<n$ previous iterates is not a satisfactory
strategy to solve nonlinear systems because it implies that $x^{k+1}$
lies in the $(p+1)$-dimensional affine subspace generated by $x^k,
x^{k-1}, \dots, x^{k-p}$. Therefore, if we repeat the
procedure~(\ref{ssmp}), we obtain that all the iterates $x^{k+1},
x^{k+2}, \dots$ lie in the same affine subspace of dimension $p+1$
related to $x^k, x^{k-1}, \dots, x^{k-p}$. Of course, this is a bad
strategy if we want to solve a nonlinear system whose solution is not
in that affine subspace.

In spite of this drawback a variation of the above formula can be used
to accelerate a pre-existent iterative process that produced the
iterates $x^{k-p}, \dots, x^k$. The idea is the following. Assume that
$x^k, \dots, x^{k-p}$ are ``previous iterates" of a given method for
solving nonlinear systems. Then, we compute, with the given method,
$x_{\trial}^{k+1}$ as a tentative candidate for being the next iterate
$x^{k+1}$. Using this trial point we compute a possible acceleration
by means of
\begin{equation} \label{posac}
x^{k+1}_{\acc} = x^k - (s^{k}_{\trial}, s^{k-1}, \dots, s^{k-p})
(y^{k}_{\trial}, y^{k-1}, \dots, y^{k-p})^{\dagger} F(x^k),
\end{equation}  
where 
\begin{equation} \label{sjyj}
s^j_{\trial} = x^{j+1}_{\trial} - x^j \mbox{ and } y^j_{\trial} =
F(x^{j+1}_{\trial}) - F(x^j) \mbox{ for all } j \leq k.
\end{equation}   
Then, if $\|F(x^{k+1}_{\acc})\|< \|F(x^{k+1}_{\trial})\|$, we define
$x^{k+1} = x^{k+1}_{\acc}$, otherwise we define $x^{k+1} =
x^{k+1}_{\trial}$. In this way, even in the case that $x^{k+1}$ is
$x^{k+1}_{\acc}$, $x^{k+1}$ does not need to lie in the affine
subspace generated by $x^k, \dots, x^{k-p}$. Therefore, the curse of
the persistent subspace is overcome.\\

\noindent
\textbf{Remark~1.} The recurrence 
\begin{equation} \label{anderson1}
x^{k+1} = x^k - (s^{k-1}_{\trial}, \dots,
s^{k-p}_{\trial})(y^{k-1}_{\trial}, \dots,
y^{k-p}_{\trial})^{\dagger} F(x^k),
\end{equation}
with $x^{j+1}_{\trial} = x^j - F(x^j)$ for all $j \leq k$ and
$s^j_{\trial}$ and $y^j_{\trial}$ defined by (\ref{sjyj}), defines
Anderson acceleration method for solving nonlinear
systems~\cite{anderson}. Many variations of this method were
introduced in recent years with theoretical justifications about the
reasons why the Anderson scheme many times accelerate the fixed point
iteration $x^{k+1} = x^k - F(x^k)$. In addition to the references
cited in the Introduction of this paper we can mention the
contributions given in
\cite{paperG,paperE,paperF,paperA,paperC,paperB,paperD}.  With the
exception of \cite{paperD}, these articles assume that the basic fixed
point iteration is convergent. In the case of \cite{paperC}, only
nonexpansiveness is assumed. In \cite{paperG}, Anderson acceleration
with local and superlinear convergence is introduced allowing
variation of the ``depth" $p$ along different iterations and
discarding stored iterates according to their relevance. In
\cite{paperE}, it is proved that Anderson acceleration improves the
convergence rate of contractive fixed-point iterations as a function
of the improvement at each step. In \cite{paperF}, the acceleration is
used to improve the Alternate Direction Multiplier Method (ADMM). In
\cite{paperA}, Anderson acceleration is interpreted as a quasi-Newton
method with GMRES-like solution of the linear systems. In \cite{paperC},
Anderson acceleration is used to improve the performance of arbitrary
minimization methods. In \cite{paperB}, it is proved that Anderson's
method is locally $r$-linearly convergent if the fixed point map is a
contraction and the coefficients in the linear combination that define
the new iteration remain bounded. In \cite{paperD}, Anderson
acceleration is applied to the solution of general nonsmooth
fixed-point problems.  Using safeguarding steps, regularization and
restarts for checking linear independence of the increments, global
convergence is proved using only nonexpansiveness of the fixed-point
iteration.\\

\noindent
\textbf{Remark~2.} The Anderson Mixing scheme is presented as an
independent method for solving nonlinear systems in~\cite[\S
  2.5]{fangsaad}. In this method one defines

\begin{equation} \label{anmix}
\bar x^{k} = x^k - (s^{k-1}, \dots, s^{k-p}) (y^{k-1}, \dots, y^{k-p})^{\dagger} F(x^k), 
\end{equation}
\begin{equation} \label{anmix2}
\bar F _k = F(x^k) - (y^{k-1}, \dots, y^{k-p}) (y^{k-1}, \dots, y^{k-p})^{\dagger} F(x^k),
\end{equation}
and
\begin{equation} \label{anmix3}
x^{k+1} = \bar x^k - \beta_k \bar F_k,
\end{equation}
where $\beta_k \equiv \beta$ for all $k$ and $\beta$ is a
problem-dependent parameter in the experiments of~\cite{fangsaad}.
So, $\bar{F}_k$ is the minimum norm element in the affine subspace
spanned by $F(x^k), \dots F(x^{k-p})$. If $p = n$ and this subspace
has dimension $n+1$ we have that $\bar F_k = 0$ and the scheme
coincides with the Sequential Secant Method.\\

The comparison between (\ref{posac}) and (\ref{anderson1}) and the
recent research on Anderson-like acceleration leads us to a natural
question: Why not to use (\ref{anderson1}) instead of (\ref{posac})
for accelerating variations of the Spectral Residual Method as
DF-SANE, defining $x^{j+1}_{\trial}$ as $x^j + \alpha_j d^j $ in
Algorithm~2.1?  The answer is that Spectral Residual methods such as
defined by Algorithm~2.1 do not possess the contractive or
non-expansive properties thar are required by all the papers that
support local and quick convergence of Anderson acceleration.

Algorithm~\ref{acceleration}.1 formally describes the way in which, at
iteration~$k$ of Algorithm~\ref{algorithm}.1 (Step~3), using an
acceleration scheme, we compute $x^{k+1}$.\\

\noindent
\textbf{Algorithm~\ref{acceleration}.1.} Let the integer number $p
\geq 1$ be given.

\begin{description}
\item[Step~1.] Define $x^{k+1}_{\trial} = x^k + \alpha_k d^k$.
\item[Step~2.] Define $\kbar = \max\{ 0, k - p + 1 \}$,
  \[
  \begin{array}{rcl}
  s^j &=& x^{j+1} - x^j \mbox{ for } j = \kbar, \dots, k-1,\\[2mm]
  y^j &=& F(x^{j+1}) - F(x^j) \mbox{ for } j = \kbar, \dots, k-1,\\[2mm]
  s^k &=& x^{k+1}_{\trial} - x^k,\\[2mm]
  y^k &=& F(x^{k+1}_{\trial}) - F(x^k),\\[2mm]
  S_k &=& (s^{\kbar}, \dots, s^{k-1}, s^k),\\[2mm]
  Y_k &=& (y^{\kbar}, \dots, y^{k-1}, y^k),
  \end{array}
  \]
  and
  \begin{equation} \label{formaqn}
    x^{k+1}_{\acc} = x^k - S_k  Y_k^\dagger F(x^k), 
  \end{equation}
  where $Y_k^\dagger$ is the Moore-Penrose pseudoinverse of $Y_k$.
\item[Step~3.] Choose $x^{k+1} \in \left\{ x^{k+1}_{\trial},
  x^{k+1}_{\acc} \right\}$ such that $\|F(x^{k+1})\| = \min \left\{
  \|F(x^{k+1}_{\trial})\|, \|F(x^{k+1}_{\acc})\| \right\}$.
\end{description}

The theorem below helps to understand the behavior of
Algorithm~\ref{algorithm}.1 with $x^{k+1}$ given by
Algorithm~\ref{acceleration}.1, called
Algorithm~\ref{algorithm}.1--\ref{acceleration}.1 from now on.
Briefly speaking, we are going to prove that, under some assumptions,
the sequence generated by
Algorithm~\ref{algorithm}.1--\ref{acceleration}.1 is such that
$\{F(x^k)\}$ converges superlinearly to zero. Thus, if~$\{x^k\}$
converges to~$x_*$ and~$J(x_*)$ is nonsingular, the convergence
of~$\{x^k\}$ to~$x_*$ is superlinear.  We do not mean that these
assumptions are ``reasonable'', in the sense that they usually, or
even frequently, hold. The theorem aims to show the correlation
between different properties of the method, which is probably useful
to understand the algorithm and, perhaps, to seek modifications and
improvements.

\begin{teo} \label{teo22}
Assume that $\{x^k\}$ is the sequence generated by
Algorithm~\ref{algorithm}.1--\ref{acceleration}.1. In addition,
suppose that
\begin{description}
\item[H1:] For all $k$ large enough we have that $x^{k+1} =
  x^{k+1}_{\acc}$.
\item[H2:] There exists a positive sequence $\beta_k \to 0$ such that
  for all $k$ large enough, the columns of $Y_k$ are linearly
  independent and
  \begin{equation} \label{difhk}
    \left\| (S_{k+1} Y_{k+1}^\dagger - S_k Y_k^\dagger) y^k \right\|
    \leq \beta_k \| y^k \|.
  \end{equation}
\item[H3:] There exists $c > 0$ such that
  \begin{equation} \label{cexiste}
    \| S_k Y_k^\dagger F(x^{k+1})\| \geq c \|F(x^{k+1})\|
  \end{equation} 
  for $k$ large enough.  
\end{description}
Then, $ \|F(x^k)\| $ converges to $0$ $Q$-superlinearly.
\end{teo}

\begin{pro}
By (\ref{formaqn}) and \textbf{H1}, we have that, for $k$ large
enough,
\begin{equation} \label{materiascura}
x^{k+1} = x^k - S_k Y_k^\dagger F(x^k).
\end{equation}
But, by the definition of $S_k$ and $Y_k$, for all $k \in \N$,
\[
S_{k+1} Y_{k+1}^\dagger y^k = s^k.
\] 
Therefore, by (\ref{difhk}) in \textbf{H2},
\begin{equation} \label{demor1}
\| S_k Y_k^\dagger y^k - s^k\| \leq \beta_k \|y^k\|.
\end{equation}
Then, by (\ref{materiascura}),
\begin{equation} \label{demor4}
\| S_k Y_k^\dagger F(x^{k+1}) \| \leq \beta_k \|y^k\|.
\end{equation}     
Thus, by (\ref{cexiste}) in \textbf{H3},
\begin{equation} \label{demor3}
c \| F(x^{k+1})\| \leq \beta_k \|y^k\|.
\end{equation}    
So, 
\[
c \|F(x^{k+1})\| \leq \beta_k (\|F(x^{k+1})\| + \|F(x^k)\|).
\]
Thus,
\[
(c - \beta_k) \|F(x^{k+1})\| \leq \beta_k \|F(x^k)\|.
\]   
Therefore, 
\[
\|F(x^{k+1})\| \leq \frac{\beta_k}{c - \beta_k} \|F(x^k)\|
\]  
for $k$ large enough, which means that $F(x^k)$ tends to zero
$Q$-superlinearly.
\end{pro}

\section{Implementation } \label{implementation}

In this section, we discuss the implementation of
Algorithm~\ref{algorithm}.1--\ref{acceleration}.1.

\subsection{ Stopping criterion}

At Step~1 of Algorithm~\ref{algorithm}.1, given $\varepsilon > 0$, we
replace the stopping criterion $\|F(x^k)\|=0$ with
\begin{equation} \label{stopcrit}
\|F(x^k)\|_2 \leq \varepsilon.
\end{equation}

\subsection{ Choice of the direction and the scaling factor}

At Step~2 of Algorithm~\ref{algorithm}.1, we must choose $\sigma_k$
and $v_k$. For this purpose, we considered the residual choice $v_k =
- F(x^k)$ for all~$k$ setting, arbitrarily, $\sigma_0=1$. A natural
choice for $\sigma_k$ ($k \geq 1$) would be the spectral step given by
\[
\sigma_k^{\spg} =
\frac{( x^k - x^{k-1} )^T ( x^k - x^{k-1} )}
     {( x^k - x^{k-1} )^T ( F(x^k) - F(x^{k-1}) )}
\]
with safeguards that guarantee that $|\sigma_k|$ belong to $[\sigmin,
  \sigmax]$.  It turns out that defining $\sigmax > 1$, preliminary
numerical experiments showed that, \textit{in the problems considered
  in the present work}, Step~2 of Algorithm~\ref{algorithm}.1 with
this choice of~$\sigma_k$ performs several backtrackings per iteration
that result in a total number of functional evaluations one order of
magnitude larger than the number of iterations; which is an unusual
behavior of nonmonotone methods based on the Barzilai-Borwein spectral
choice~\cite{bmr,bmr2,bmrsurveyspg,bmr5,lmr,lacruzraydan,raydan1,raydan2}. On
the other hand, it was also observed that, at Step~3 of
Algorithm~\ref{acceleration}.1, the acceleration step was always
chosen.  This means that the costly obtained $x^{k+1}_{\trial}$ was
always beaten by $x^{k+1}_{\acc}$. These two observations suggested
that a more conservative, i.e.\ small, scaling factor~$\sigma_k$
should be considered. A small $\sigma_k$ could result in a trial point
$x^{k+1}_{\trial} = x^k + \alpha_k d^k$ with $\alpha_k=\pm 1$ that
satisfies~(\ref{algarmijo}), i.e.\ no backtracking, and that provides
the information required to compute a successful $x^{k+1}_{\acc}$ at
Step~2 of Algorithm~\ref{acceleration}.1.

The conservative choice of $\sigma_k$ employed in our implementation
is as follows. We consider an algorithmic parameter
\begin{equation} \label{hinit}
h_{\mathrm{init}} > 0
\end{equation}
and we define:
\begin{equation}\label{sigmima}
\sigmin = \sqrt{\epsilon}, \; \sigmax=1, 
\end{equation} 
\begin{equation} \label{sigma0}
\sigma_0 = 1,
\end{equation}
\begin{equation} \label{epslambda1}
\bar \sigma_k = h_{\mathrm{init}} \frac{\| x^k - x^{k-1} \|}{\| F(x^k)\|} \mbox{ for all } k \geq 1, 
\end{equation} 
and
\begin{equation} \label{sigmundo}
\sigma_k = \bar \sigma_k \mbox{ if } k \geq 1 \mbox{ and } \bar
\sigma_k \in [\max\{1,\|x^k\|\} \sigmin,\sigmax].
\end{equation}
Finally, if $k \geq 1$ and $\bar \sigma_k \notin [\max\{1,\|x^k\|\}
  \sigmin,\sigmax]$, we define
 \begin{equation} \label{epslambda2}
\bar{\bar \sigma}_k = h_{\mathrm{init}} \frac{\| x^k \|}{\| F(x^k)\|}
\end{equation}    
and we compute $\sigma_k$ as the projection of $ \bar{\bar \sigma}_k$
onto the interval $[\max\{1,\|x^k\|\} \sigmin,\sigmax]$.


\subsection{Procedure for reducing the step}

At Step~2.4 of Algorithm~\ref{algorithm}.1, we compute
$\alpha_+^{\new}$ as the minimizer of the the univariate quadratic
$q(\alpha)$ that interpolates $q(0)=f(x^k)$,
$q(\alpha_+)=f(x^k-\alpha_+ \sigma_k F(x^k))$, and $q'(0)=-\sigma_k
F(x^k)^T \nabla f(x^k) = -\sigma_k F(x^k)^T J(x^k)
F(x^k)$. Following~\cite{lmr}, since we consider $J(x^k)$ unavailable,
we consider $J(x^k) = I$. Thus,
\[
\alpha_+^{\new} = \max \left\{ \taumin \alpha_+, \min \left\{
\frac{\alpha_+^2 f(x^k)}{f(x^k - \alpha_+ \sigma_k F(x^k)) +
  (2\alpha_+-1) f(x^k)}, \taumax \alpha_+ \right\} \right\}.
\]
Analogously,
\[
\alpha_-^{\new} = \max \left\{ \taumin \alpha_-, \min \left\{
\frac{\alpha_-^2 f(x^k)}{f(x^k + \alpha_- \sigma_k F(x^k)) +
  (2\alpha_--1) f(x^k)}, \taumax \alpha_- \right\} \right\}.
\]

\subsection{Computing the acceleration} \label{secQR}

In Step~2 of Algorithm~\ref{acceleration}.1, computing
$x^{k+1}_{\acc}$ as defined in~(\ref{formaqn}) is equivalent to
computing the minimum-norm least-squares solution~$\bar \omega$ to the
linear system $Y_k \omega = F(x^{k+1}_{\trial})$ and defining
$x^{k+1}_{\acc} = x^{k+1}_{\trial} - S_k \bar \omega$. In practice, we
compute the minimum-norm least-squares solution with a complete
orthogonalization of $Y_k$. Computing this factorization from scratch
at each iteration could be expensive. However, by definition, $Y_k$
corresponds to a slight variation of $Y_{k-1}$. In practice, when
$Y_k$ has not full numerical rank, one extra column is added to it;
and if $Y_k$ has null numerical rank, a new $Y_k$ is computed from
scratch. For completeness, the practical implementation of
Algorithm~\ref{acceleration}.1 is given below.\\

\noindent
\textbf{Algorithm~\ref{implementation}.1.} Let
$0<h_{\mathrm{small}}<h_{\mathrm{large}}$ and $p \geq 1$ be given. Set
$r_{\max} \leftarrow 0$ and $\ell \leftarrow 1$.

\begin{description}
\item[Step~1.] Define $x^{k+1}_{\trial} = x^k + \alpha_k d^k$.
\item[Step~2.] 
\item[Step~2.1.] If $k=0$, then set $S_k$ and $Y_k$ as
  matrices with zero columns. Otherwise, set $S_k \leftarrow S_{k-1}$
  and $Y_k \leftarrow Y_{k-1}$.
\item[Step~2.2.] If $S_k$ and $Y_k$ have $p$ columns,
  then remove the leftmost column of $S_k$ and $Y_k$.
\item[Step~2.3.] Add $s^k = x^{k+1}_{\trial} - x^k$ and $y^k =
  F(x^{k+1}_{\trial}) - F(x^k)$ as rightmost column of $S_k$ and $Y_k$,
  respectively; and set $r_{\max} \leftarrow \max\{ r_{\max},
  \mathrm{rank}(Y_k) \}$.
\item[Step~2.4.] If $\mathrm{rank}(Y_k) < r_{\max}$, then
  execute Steps~2.4.1--2.4.2.
  \begin{description}
  \item[Step~2.4.1.] If $S_k$ and $Y_k$ have $p$
    columns, then remove the leftmost column of $S_k$ and
    $Y_k$.
  \item[Step~2.4.2.] Add $\sextra = \xextra - x^k$ and $\yextra =
    F(\xextra) - F(x^k)$ as rightmost column of $S_k$ and $Y_k$,
    respectively, where $\xextra = x^k + h_{\mathrm{small}} \,
    e_\ell$. Set $r_{\max} \leftarrow \max\{ r_{\max},
    \mathrm{rank}(Y_k) \}$ and $\ell \leftarrow
    \mathrm{mod}(\ell,n)+1$.
  \end{description}
\item[Step~2.5.] If $\mathrm{rank}(Y_k) \neq 0$, then execute
  Steps~2.5.1--2.5.3.
  \begin{description}
    \item[Step~2.5.1.] Compute the minimum-norm least-squares solution
      $\bar \omega$ to the linear system $Y_k \omega =
      F(x^k)$ and define $x^{k+1}_{\acc} =
      x^k - S_k \bar \omega$.
    \item[Step~2.5.2.] If Step~2.4.2 was executed, then
      remove the rightmost column of $S_k$ and $Y_k$, i.e.\ columns
      $\sextra$ and $\yextra$.
    \item[Step~2.5.3.] If $x^{k+1}_{\acc} \neq x^k$,
      $\|x^{k+1}_{\acc}\| \leq 10 \max\{1,\|x^k\|\}$, and $\|
      F(x^{k+1}_{\acc}) \| < \| F(x^{k+1}_{\trial}) \|$, then redefine
      $x^{k+1}_{\trial} = x^{k+1}_{\acc}$, substitute the rightmost
      column of $S_k$ and $Y_k$, i.e.\ columns $s^k$ and $y^k$, with
      $\sacc = x^{k+1}_{\acc} - x^k$ and $\yacc = F(x^{k+1}_{\acc}) -
      F(x^k)$, respectively, and set $r_{\max} \leftarrow \max\{
      r_{\max}, \mathrm{rank}(Y_k) \}$.
  \end{description}
\item[Step~2.6.] If $\mathrm{rank}(Y_k) = 0$, then execute
  Steps~2.6.1--2.6.3.
  \begin{description}
    \item[Step~2.6.1.] Redefine matrices $S_k$ and $Y_k$ as
      matrices with zero columns.
    \item[Step~2.6.2.] Execute Steps~2.6.2.1 $p-1$ times.
      \begin{description}
      \item[Step~2.6.2.1.] Add $\sextra = \xextra - x^{k+1}_{\trial}$
        and $\yextra = F(\xextra) - F(x^{k+1}_{\trial})$ as rightmost
        column of $S_k$ and $Y_k$, respectively, where $\xextra = x^k
        + h_{\mathrm{large}} \, e_\ell$. Set $r_{\max} \leftarrow
        \max\{ r_{\max}, \mathrm{rank}(Y_k) \}$ and $\ell \leftarrow
        \mathrm{mod}(\ell,n)+1$.
      \end{description}
    \item[Step~2.6.3.] Add $s^k = x^{k+1}_{\trial} - x^k$ and
      $y^k = F(x^{k+1}_{\trial}) - F(x^k)$ as rightmost column of $S_k$ and
      $Y_k$, respectively; and set $r_{\max} \leftarrow \max\{
      r_{\max}, \mathrm{rank}(Y_k) \}$.
  \end{description}
\item[Step~2.7.] Compute the minimum-norm least-squares
  solution $\bar \omega$ to the linear system $Y_k \omega = F(x^k)$ and define
  $x^{k+1}_{\acc} = x^k - S_k \bar \omega$.
\item[Step~2.8.] If $x^{k+1}_{\acc} \neq x^k$, $\|x^{k+1}_{\acc}\|
  \leq 10 \max\{1,\|x^k\|\}$, and $\| F(x^{k+1}_{\acc}) \| < \|
  F(x^{k+1}_{\trial}) \|$, then redefine $x^{k+1}_{\trial} =
  x^{k+1}_{\acc}$ and substitute the rightmost column of $S_k$ and
  $Y_k$, i.e.\ columns $s^k$ and $y^k$, with $\sacc = x^{k+1}_{\acc} -
  x^k$ and $\yacc = F(x^{k+1}_{\acc}) - F(x^k)$, respectively, and set
  $r_{\max} \leftarrow \max\{ r_{\max}, \mathrm{rank}(Y_k) \}$.
\item[Step~3.] Define $x^{k+1} = x^{k+1}_{\trial}$.
\end{description}

In general, Step~2.6 is not executed. If iteration $k$ is such that
Step~2.6 was not executed in the previous~$p$ iterations, then
matrices $S_k$ and $Y_k$ correspond to removing the leftmost column
from and adding a new rightmost column to $S_{k-1}$ and $Y_{k-1}$,
respectively. The leftmost column is not removed if the maximum number
of columns $p$ was not reached yet. When $k=0$, $Y_k$ is a
single-column matrix for which matrices $P$, $Q$, and~$R$ of the
rank-revealing QR decomposition $Y_k P = Q R$ can be trivially
computed. For $k>0$, the QR decomposition of $Y_k$ can be obtained
with time complexity $O(n \, p)$ by updating, via Givens rotations,
the QR decomposition of $Y_{k-1}$. Moreover, by using the QR
decomposition of $Y_k$, the minimum-norm least-squares solution $\bar
\omega$ can be computed with time complexity $O(n + p^2)$ if $Y_k$ is
full-rank and with time complexity $O(n + p^3)$ in the rank-deficient
case. Summing up, by assuming that $p \ll n$ and that $p$ does not
depend on $n$, iterations of Algorithm~\ref{implementation}.1 can be
implemented with time complexity $O(n)$. The space complexity is $O(n
\, p + p^2)$ and it is related to the fact that we must save matrix
$S_k \in \R^{n \times p}$ and matrices $Q \in \R^{n \times p}$ and $R
\in \R^{p \times p}$ of the QR decomposition of $Y_k$. Of course, the
permutation matrix $P$ can be saved in an array of size $p$. Thus, the
space complexity of Algorithm~\ref{implementation}.1 is also $O(n)$
under the same assumptions. When Step~2.6 is executed, the QR
decomposition of $Y_k$ must be computed from scratch, with time
complexity~$O(n \, p^2)$.

\section{Numerical experiments} \label{experiments}

We implemented Algorithm~\ref{algorithm}.1 and
Algorithm~\ref{implementation}.1 (the practical version of
Algorithm~\ref{acceleration}.1) in Fortran~90. In the numerical
experiments, we considered the standard values $\gamma=10^{-4}$,
$\taumin=0.1$, $\taumax=0.5$, $M=10$, $\sigmin=\sqrt{\epsilon}$,
$\sigmax=1/\sqrt{\epsilon}$, and $\eta_k = 2^{-k} \min\{ \frac{1}{2}
\| F(x^0)\|, \sqrt{\| F(x^0) \| }\}$, where $\epsilon \approx
10^{-16}$ is the machine precision. More crucial to the method
performance are the choices of $h_{\mathrm{init}}$,
$h_{\mathrm{small}}$, $h_{\mathrm{large}}$, and $p$, whose values are
mentioned below. All tests were conducted on a computer with a 3.4 GHz
Intel Core i5 processor and 8GB 1600 MHz DDR3 RAM memory, running
macOS Mojave (version 10.14.6). Code was compiled by the GFortran
compiler of GCC (version 8.2.0) with the -O3 optimization directive
enabled.

\subsection{2D and 3D Bratu problem} \label{sec61}

In a first experiment, we considered 2D and 3D versions of the Bratu
problem
\begin{equation} \label{bratu}
- \Delta u + \theta e^u = \phi(\bar u) \mbox{ in } \Omega
\end{equation}
with boundary conditions $u=\bar u$ on $\partial \Omega$. In the 2D
case, $\Omega=[0,1]^2$ and, following~\cite{kelley95}, we set $\bar u
= 10 u_1 u_2 (1-u_1) (1-u_2) e^{u_1^{4.5}}$; while in the 3D case,
$\Omega=[0,1]^3$ and $\bar u = 10 u_1 u_2 u_3 (1-u_1) (1-u_2) (1-u_3)
e^{u_1^{4.5}}$. In both cases, $\phi(u)=-\Delta u + \theta e^u$; so
the problem has $\bar u $ as known solution. Considering $n_p$
discretization points in each dimension and approximating
\[
\Delta u(x) \approx \frac{u(x \pm h e_1) + u(x \pm h e_2) - 4 u(x)}{h^2} 
\]
and
\[
\Delta u(x) \approx \frac{u(x \pm h e_1) + u(x \pm h e_2) + u(x \pm h e_3) - 6 u(x)}{h^2},
\]
where $h=1/(n_p-1)$ and $e_i$ is the $i$th canonical vector in the
corresponding space ($\R^2$ or $\R^3$), we obtain nonlinear systems of
equations with $n=(n_p-2)^2$ and $n=(n_p-2)^3$ variables in the 2D and
3D cases, respectively. Starting from $u=0$, fixing $\theta = -100$,
and varying $n_p \in \{ 100, 125, \dots, 400 \}$ and $n_p \in \{ 10,
15, \dots, 70 \}$ in the 2D and 3D cases respectively, we run NITSOL
(file \textsc{nitsol,11-1-05.tar.gz} downloaded from
\url{https://users.wpi.edu/~walker/NITSOL/} on October 26th, 2020) and
the method proposed in the present work, which will be called
Accelerated DF-SANE from now on. As stopping criterion, we
considered~(\ref{stopcrit}) with $\varepsilon=10^{-6} \sqrt{n}$. For
NITSOL, we use all its default parameters, except the maximum number
of (nonlinear) iterations, which was increased in order to avoid
premature stops. By default, NITSOL corresponds to an Inexact-Newton
method in which Newtonian systems are solved with GMRES (maximum
Krylov subspace dimension equal to~$20$), approximating the
Jacobian-vector products by finite-differences. For Accelerated
DF-SANE, based on preliminary experiments, we set
$h_{\mathrm{small}}=10^{-4}$, $h_{\mathrm{large}}=0.1$, and
$h_{\mathrm{init}}=0.01$ in the 2D case and $h_{\mathrm{small}} =
h_{\mathrm{large}}=0.1$ and $h_{\mathrm{init}}=1$ in the 3D case. In
both cases, we considered~$p=5$.

Tables~\ref{tab1} and~\ref{tab2} show the results. In the tables,
\#it$_1$ corresponds to the nonlinear iterations (outer iterations) of
NITSOL; while \#it$_2$ corresponds to its linear iterations (inner or
GMRES iterations). For both methods, ``fcnt'' stands for number of
evaluations of $F$, ``Time'' stands for CPU time in seconds, and
``SC'' stands for ``Stopping criterion''. Remaining columns are
self-explanatory. In the case of NITSOL, stopping criterion equal
to~$0$ means success; while~$6$ means ``failure to reach an acceptable
step through backtracking''. Accelerated DF-SANE satisfied the
stopping criterion related to success in all problems. Regarding the
2D problems, it should be first noted that NITSOL failed in satisfying
the stopping criterion for $n_p \geq 225$. Accelerated DF-SANE is
faster than NITSOL in all problems that both methods solved. In the
larger problem that both methods solved, Accelerated DF-SANE is around
thirty times faster than NITSOL. Regarding the 3D problems
(Table~\ref{tab2}), both methods satisfied the stopping criterion
related to success in all problems.  Accelerated DF-SANE is faster
than NITSOL in all problems in the table. In the larger problem in the
table ($n_p=70$), Accelerated DF-SANE is more than twenty times faster
than NITSOL.

\begin{table}[ht!]
\begin{center}
\resizebox{0.85\textwidth}{!}{
\begin{tabular}{|cr|ccrrrr|crrr|}
\hline
\multirow{2}{*}{$n_p$} & \multicolumn{1}{c|}{\multirow{2}{*}{$n$}} &
\multicolumn{6}{c}{NITSOL (Newton-GMRES)} & \multicolumn{4}{|c|}{Accelerated DF-SANE} \\
\cline{3-12}
& &
\multicolumn{1}{|c}{SC} & \multicolumn{1}{c}{$\|F(x_*)\|_2$} &
\multicolumn{1}{c}{\#it$_1$} & \multicolumn{1}{c}{\#it$_{2}$} & \multicolumn{1}{c}{fcnt} &
\multicolumn{1}{c}{Time} & \multicolumn{1}{|c}{$\|F(x_*)\|_2$} & \multicolumn{1}{c}{\#it} &
\multicolumn{1}{c}{fcnt} & \multicolumn{1}{c|}{Time} \\
\hline
\hline
100 &    9,604 & 0 & 9.7e$-$05 &    220 &  197,732 &  197,953 &   39.63 & 9.8e$-$05 &  5,219 &  10,688 &    4.87 \\  
125 &   15,129 & 0 & 1.2e$-$04 &    631 &  473,616 &  474,248 &  151.88 & 1.2e$-$04 &  2,612 &   5,489 &    3.80 \\  
150 &   21,904 & 0 & 1.5e$-$04 &    379 &  315,338 &  315,718 &  165.07 & 1.5e$-$04 &  2,946 &   6,007 &    6.59 \\  
175 &   29,929 & 0 & 1.7e$-$04 &    401 &  331,380 &  331,782 &  265.53 & 1.7e$-$04 &  4,475 &  10,007 &   17.46 \\  
200 &   39,204 & 0 & 2.0e$-$04 &  4,421 &  846,957 &  851,385 &  983.97 & 2.0e$-$04 &  6,775 &  14,385 &   32.24 \\  
225 &   49,729 & 6 & 3.3e$+$01 &  3,944 &  173,310 &  177,283 &  261.03 & 2.2e$-$04 &  4,328 &   8,927 &   26.72 \\  
250 &   61,504 & 6 & 4.1e$+$01 &    105 &   46,897 &   47,019 &   82.69 & 2.5e$-$04 & 12,661 &  26,353 &   94.47 \\  
275 &   74,529 & 6 & 5.3e$+$01 &  1,411 &   94,194 &   95,635 &  224.74 & 2.7e$-$04 &  8,809 &  19,583 &   88.07 \\  
300 &   88,804 & 6 & 7.2e$+$01 &  1,341 &  104,735 &  106,112 &  317.71 & 3.0e$-$04 & 15,858 &  34,194 &  190.49 \\  
325 &  104,329 & 6 & 9.6e$+$01 &  2,430 &  194,640 &  197,126 &  705.53 & 3.2e$-$04 & 10,791 &  23,403 &  160.02 \\  
350 &  121,104 & 6 & 1.4e$+$02 &    195 &   82,940 &   83,151 &  319.00 & 3.5e$-$04 & 10,805 &  25,915 &  176.52 \\  
375 &  139,129 & 6 & 1.7e$+$02 &  1,127 &  121,100 &  122,243 &  649.94 & 3.7e$-$04 & 16,335 &  38,648 &  331.03 \\  
400 &  158,404 & 6 & 2.2e$+$02 &  2,936 &  144,460 &  147,415 &  851.67 & 4.0e$-$04 & 21,106 &  55,901 &  515.65 \\  
\hline
\end{tabular}}
\end{center}
\caption{Performances of NITSOL and Accelerated DF-SANE in the 2D Bratu problem.}
\label{tab1}
\end{table}

\begin{table}[ht!]
\begin{center}
\resizebox{0.85\textwidth}{!}{
\begin{tabular}{|cr|crrrr|crrr|}
\hline
\multirow{2}{*}{$n_p$} & \multicolumn{1}{c|}{\multirow{2}{*}{$n$}} &
\multicolumn{5}{c}{NITSOL (Newton-GMRES)} & \multicolumn{4}{|c|}{Accelerated DF-SANE} \\
\cline{3-11}
& &
\multicolumn{1}{c}{$\|F(x_*)\|_2$} &
\multicolumn{1}{c}{\#it$_1$} & \multicolumn{1}{c}{\#it$_{2}$} & \multicolumn{1}{c}{fcnt} &
\multicolumn{1}{c}{Time} & \multicolumn{1}{|c}{$\|F(x_*)\|_2$} & \multicolumn{1}{c}{\#it} &
\multicolumn{1}{c}{fcnt} & \multicolumn{1}{c|}{Time} \\
\hline
\hline
 10 &      512 & 1.7e$-$05 &     5 &      213 &      219 &    0.00 & 2.1e$-$05 &    126 &     308 &    0.00 \\  
 15 &    2,197 & 3.7e$-$05 &     7 &    1,621 &    1,629 &    0.07 & 4.7e$-$05 &    223 &     662 &    0.05 \\  
 20 &    5,832 & 6.1e$-$05 &    11 &    7,157 &    7,169 &    0.97 & 7.2e$-$05 &  1,221 &   4,271 &    0.92 \\  
 25 &   12,167 & 8.8e$-$05 &    19 &   15,893 &   15,913 &    4.62 & 1.1e$-$04 &    529 &   1,840 &    0.78 \\  
 30 &   21,952 & 1.2e$-$04 &    25 &   21,991 &   22,017 &   12.72 & 1.4e$-$04 &    812 &   3,012 &    2.29 \\  
 35 &   35,937 & 1.6e$-$04 &    38 &   34,935 &   34,974 &   36.31 & 1.9e$-$04 &  1,210 &   4,530 &    6.37 \\  
 40 &   54,872 & 2.2e$-$04 &    70 &   66,115 &   66,186 &  118.56 & 2.3e$-$04 &  1,109 &   4,379 &    9.35 \\  
 45 &   79,507 & 2.8e$-$04 &   147 &  137,626 &  137,774 &  386.09 & 2.8e$-$04 &  1,315 &   5,444 &   20.06 \\  
 50 &  110,592 & 3.3e$-$04 &   217 &  199,042 &  199,260 &  759.61 & 3.3e$-$04 &  1,574 &   6,501 &   30.38 \\  
 55 &  148,877 & 3.8e$-$04 &   357 &  312,023 &  312,381 & 1736.95 & 3.8e$-$04 &  1,706 &   7,254 &   53.10 \\  
 60 &  195,112 & 4.4e$-$04 &   504 &  418,201 &  418,706 & 3074.25 & 4.4e$-$04 &  1,821 &   8,019 &   70.29 \\  
 65 &  250,047 & 5.0e$-$04 &   988 &  691,192 &  692,181 & 6715.18 & 4.9e$-$04 &  2,061 &   9,379 &  109.39 \\  
 70 &  314,432 & 5.6e$-$04 &   609 &  486,651 &  487,261 & 5735.90 & 5.6e$-$04 &  1,836 &   8,431 &  118.13 \\  
\hline
\end{tabular}}
\end{center}
\caption{Performances of NITSOL and Accelerated DF-SANE in the 3D Bratu problem.}
\label{tab2}
\end{table}

\begin{figure}[ht!]
\begin{center}
\resizebox{0.5\textwidth}{!}{\input{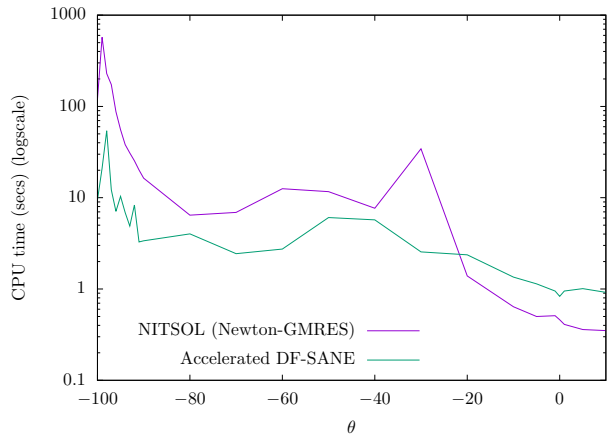}}
\vspace{-0.8cm}
\end{center}
\caption{This figure shows the CPU time (in seconds) used by NITSOL
  and Accelerated DF-SANE to solve the 3D Bratu problem with $n_p=40$ and
  $\theta \in [-100,10]$.}
\label{fig1}
\end{figure}

The difficulty of the Bratu problem varies with the value of $\theta$,
that may be positive or negative. For the formulation given
in~(\ref{bratu}), negative values of $\theta$ correspond to more
difficult problem. Therefore, results in Tables~\ref{tab1}
and~\ref{tab2} raise the question of how the performances of the
methods compare to each other for different values of~$\theta$. So, we
arbitrarily considered the 3D Bratu problem with $n_p=40$, which is
affordable for both methods and varied $\theta \in
[-100,10]$. Figure~\ref{fig1} shows the results of applying NITSOL and
Accelerated DF-SANE. The graphic shows that for $\theta \geq -20$ both
methods use less than a second and NITSOL outperforms Accelerated
DF-SANE. On the other hand, in the most difficult problems
(i.e.\ $\theta < -20$), where up to one thousand seconds are required,
Accelerated DF-SANE outperforms NITSOL by approximately an order of
magnitude. (Note that the $y$-axis is in logarithmic scale.)
Considering the number of functional evaluations as a performance
metric, instead of the CPU time, results are mostly the same.

\begin{figure}[ht!]
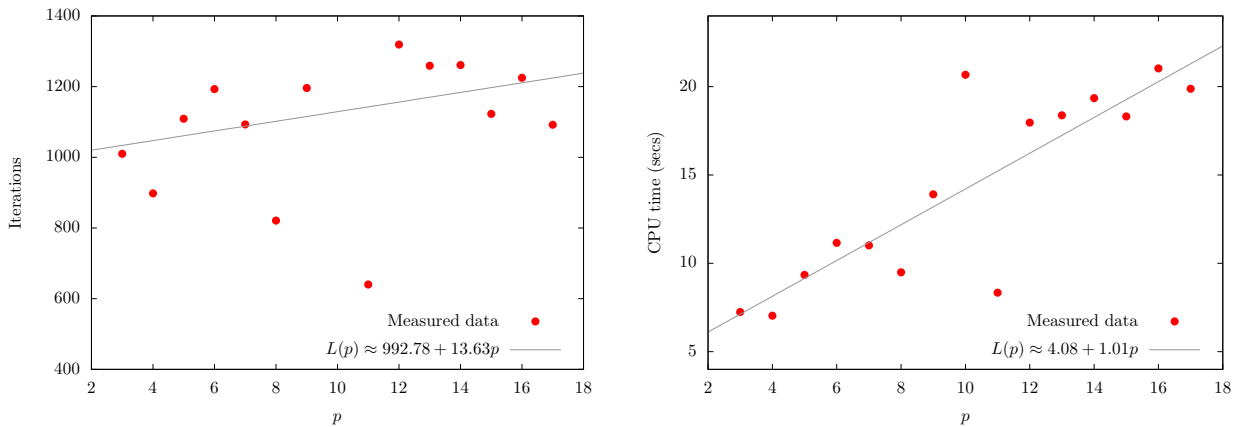

\begin{center}
\resizebox{\textwidth}{!}{\begin{tabular}{cc}
\input{bmdiisfig2a.tex} & \input{bmdiisfig2b.tex}
\end{tabular}}
\vspace{-0.8cm}
\end{center}
\caption{Performance of Accelerated DF-SANE in the 3D Bratu problem with
  $n_p=40$ and $\theta=-100$ varying the number $p$ of columns in
  matrices $S_k$ and $Y_k$.}
\label{fig2}
\end{figure}

Another relevant question is how the performance of Accelerated
DF-SANE varies as a function of its parameter~$p$. So, considering
again the 3D Bratu problem with $\theta=-100$ and $n_p=40$, we ran
Accelerated DF-SANE with $p \in \{ 3, 4, \dots, 17
\}$. Figure~\ref{fig2} shows the results. The figure on the left shows
that the number of iterations of Accelerated DF-SANE is nearly
constant as a function of~$p$; while the figure on the right shows, as
expected, that, the larger~$p$, the larger the CPU time per
iteration. Note that for the considered values of $p$, the variation
on the number of iterations is up to 20\%; while the variation in CPU
time is up to 300\%.

\begin{figure}[ht!]
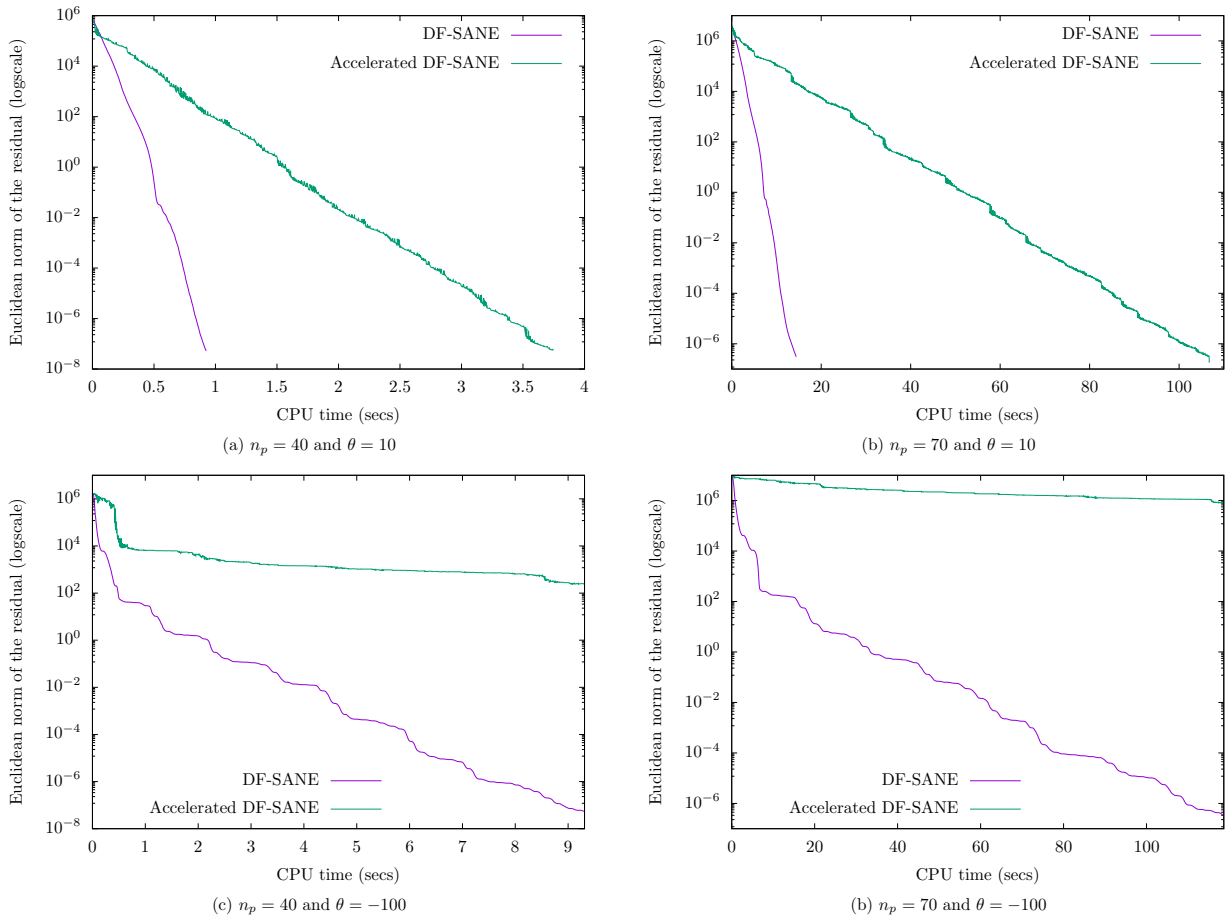

\begin{center}
\resizebox{\textwidth}{!}{\begin{tabular}{cc}
\input{bmdiisfig31ab.tex} & \input{bmdiisfig32ab.tex} \\
(a) $n_p=40$ and $\theta=10$ & (b) $n_p=70$ and $\theta=10$ \\
\input{bmdiisfig33ab.tex} & \input{bmdiisfig34ab.tex} \\
(c) $n_p=40$ and $\theta=-100$ & (b) $n_p=70$ and $\theta=-100$ \\
\end{tabular}}
\end{center}
\caption{Performance of DF-SANE and Accelerated DF-SANE in four instances of
  the 3D Bratu problem with $\theta \in \{-100,10\}$ and $n_p \in
  \{40,70\}$.}
\label{fig3}
\end{figure}

Another relevant question regarding Accelerated DF-SANE is how it
compares against DF-SANE, i.e.\ the original method without the
acceleration being proposed in the present work. We were unable to run
DF-SANE in the set of problems considered in Tables~\ref{tab1}
and~\ref{tab2}, because DF-SANE is unable to reach the stopping
criterion~(\ref{stopcrit}) with $\varepsilon=10^{-6} \sqrt{n}$ within
an affordable time. As an alternative, we considered the four
instances of the 3D Bratu problem with $n_p \in \{40,70\}$ and $\theta
\in \{-100,10\}$. In the two instances with $\theta=10$, DF-SANE was
able to satisfy the imposed stopping criterion within an affordable
time. This result was in fact expected because for this value of
$\theta$ the Bratu problem resembles the minimization of a convex
quadratic function; and Barzilai-Borwein or spectral step based
methods are expected to perform especially well in this case. In the
two instances with $\theta=-100$, we first run Accelerated DF-SANE and
then we run DF-SANE using as CPU time limit the time consumed by
Accelerated DF-SANE. Figure~\ref{fig3} shows the results. It is very
clear from the four graphics that the acceleration process is very
efficient in its purpose of accelerating DF-SANE. (In the cases with
$\theta=-100$, with an excruciatingly slow progress, DF-SANE is far
from reaching convergence after a couple of hours of CPU time.)

\subsection{Modified Bratu, Driven Cavity, and Flow in a Porous Media problems} \label{sec62}

In a second set of experiments, we considered three PDE-based problems
described in~\cite{pernice}. Fortran implementations of these problems
are included as examples of usage in the NITSOL distribution. The
three problems are discretizations of 2D PDE-based systems of nonlinear
equations. All parameters of the problems were set as suggested
in~\cite{pernice}. We varied the value of $n_p$ for the three problems
trying to illustrate the behavior of the methods and solving as large
as possible problems with both methods within an affordable
time. NITSOL was run with all its default parameters, as in the
previous section; while Accelerated DF-SANE was run with the same
parameters already reported for the 2D Bratu problem in the previous
section. Tables~\ref{tab3}--\ref{tab5} show the
results. Table~\ref{tab3} shows that NITSOL failed in satisfying the
stopping criterion for $n_p \geq 135$ in the Driven Cavity
problem. Accelerated DF-SANE was faster than NITSOL in all instances
that both methods solved. In the largest instance that both methods
solved, Accelerated DF-SANE is more than thirty times faster than
NITSOL. Table~\ref{tab4} shows that NITSOL is faster than Accelerated
DF-SANE in the smaller instances of the Generalized 2D Bratu problem
($n_p$ up to $750$). In the largest instance in the table
($n_p=1,500$), Accelerated DF-SANE takes around half of the time
required by NITSOL. Finally, Table~\ref{tab5} shows that NITSOL is
faster than Accelerated DF-SANE in the smaller instances ($n_p$ up to
$400$) of the Flow in a Porous Media problem. In the largest instance
in the table ($n_p=1,500$), Accelerated DF-SANE takes around half of
the time required by NITSOL.

\begin{table}[ht!]
\begin{center}
\resizebox{0.85\textwidth}{!}{
\begin{tabular}{|cr|ccrrrr|crrr|}
\hline
\multirow{2}{*}{$n_p$} & \multicolumn{1}{c|}{\multirow{2}{*}{$n$}} &
\multicolumn{6}{c}{NITSOL (Newton-GMRES)} & \multicolumn{4}{|c|}{Accelerated DF-SANE} \\
\cline{3-12}
& &
\multicolumn{1}{|c}{SC} & \multicolumn{1}{c}{$\|F(x_*)\|_2$} &
\multicolumn{1}{c}{\#it$_1$} & \multicolumn{1}{c}{\#it$_{2}$} & \multicolumn{1}{c}{fcnt} &
\multicolumn{1}{c}{Time} & \multicolumn{1}{|c}{$\|F(x_*)\|_2$} & \multicolumn{1}{c}{\#it} &
\multicolumn{1}{c}{fcnt} & \multicolumn{1}{c|}{Time} \\
\hline
\hline
  10 &      100 & 0 & 7.7e$-$06 &      5 &      194 &       200 &    0.00 & 9.6e$-$06 &     96 &     200 &    0.00 \\         
  15 &      225 & 0 & 1.0e$-$05 &      6 &      988 &       995 &    0.00 & 1.4e$-$05 &    190 &     392 &    0.00 \\         
  20 &      400 & 0 & 1.6e$-$05 &      6 &    2,739 &     2,746 &    0.02 & 1.9e$-$05 &    278 &     582 &    0.01 \\         
  25 &      625 & 0 & 2.0e$-$05 &     10 &    6,273 &     6,284 &    0.09 & 2.4e$-$05 &    473 &   1,044 &    0.03 \\         
  30 &      900 & 0 & 2.4e$-$05 &     15 &   11,939 &    11,955 &    0.24 & 2.9e$-$05 &    663 &   1,519 &    0.06 \\         
  35 &    1,225 & 0 & 3.1e$-$05 &     23 &   20,058 &    20,082 &    0.57 & 3.4e$-$05 &    889 &   2,118 &    0.11 \\         
  40 &    1,600 & 0 & 3.8e$-$05 &     36 &   32,819 &    32,856 &    1.19 & 3.9e$-$05 &  1,141 &   2,954 &    0.20 \\         
  45 &    2,025 & 0 & 4.4e$-$05 &     56 &   52,200 &    52,257 &    2.52 & 4.4e$-$05 &  1,431 &   3,706 &    0.33 \\         
  50 &    2,500 & 0 & 5.0e$-$05 &     80 &   74,981 &    75,062 &    4.68 & 5.0e$-$05 &  1,747 &   4,653 &    0.50 \\         
  55 &    3,025 & 0 & 5.2e$-$05 &    114 &  106,481 &   106,596 &    8.27 & 5.5e$-$05 &  2,093 &   6,010 &    0.75 \\         
  60 &    3,600 & 0 & 5.8e$-$05 &    154 &  142,222 &   142,377 &   13.09 & 6.0e$-$05 &  2,470 &   7,346 &    1.06 \\         
  65 &    4,225 & 0 & 6.3e$-$05 &    217 &  195,849 &   196,067 &   21.29 & 6.5e$-$05 &  2,880 &   8,824 &    1.48 \\         
  70 &    4,900 & 0 & 6.8e$-$05 &    278 &  244,149 &   244,428 &   30.60 & 7.0e$-$05 &  3,297 &  10,533 &    1.99 \\         
  75 &    5,625 & 0 & 7.4e$-$05 &    392 &  328,630 &   329,023 &   48.01 & 7.5e$-$05 &  3,781 &  12,557 &    2.74 \\         
  80 &    6,400 & 0 & 7.8e$-$05 &    568 &  442,630 &   443,199 &   72.66 & 8.0e$-$05 &  4,274 &  14,564 &    3.48 \\         
  85 &    7,225 & 0 & 8.4e$-$05 &    790 &  563,030 &   563,821 &  104.84 & 8.5e$-$05 &  4,811 &  16,905 &    4.47 \\         
  90 &    8,100 & 0 & 9.0e$-$05 &  1,135 &  706,031 &   707,167 &  147.14 & 9.0e$-$05 &  5,369 &  19,423 &    5.67 \\         
  95 &    9,025 & 0 & 9.5e$-$05 &  1,661 &  856,631 &   858,293 &  200.43 & 9.5e$-$05 &  5,949 &  22,052 &    7.10 \\         
 100 &   10,000 & 0 & 1.0e$-$04 &  3,051 & 1040,211 & 1,043,263 &  269.07 & 1.0e$-$04 &  6,563 &  24,901 &    8.76 \\         
 105 &   11,025 & 0 & 1.0e$-$04 & 13,222 & 1233,251 & 1,246,474 &  357.23 & 1.0e$-$04 &  7,208 &  27,977 &   10.74 \\         
 110 &   12,100 & 0 & 1.1e$-$04 & 20,017 & 1415,812 & 1,435,830 &  452.47 & 1.1e$-$04 &  7,593 &  30,394 &   12.62 \\         
 115 &   13,225 & 0 & 1.1e$-$04 & 18,717 & 1290,772 & 1,309,495 &  458.66 & 1.1e$-$04 &  8,543 &  34,375 &   15.81 \\         
 120 &   14,400 & 0 & 1.2e$-$04 & 18,350 & 1461,992 & 1,480,349 &  568.50 & 1.2e$-$04 &  9,265 &  37,904 &   18.57 \\         
 125 &   15,625 & 0 & 1.2e$-$04 & 26,341 & 1628,992 & 1,655,343 &  696.52 & 1.2e$-$04 & 10,036 &  41,724 &   22.07 \\         
 130 &   16,900 & 0 & 1.3e$-$04 & 13,114 & 1895,392 & 1,908,981 &  903.98 & 1.3e$-$04 & 10,824 &  45,718 &   26.53 \\         
 135 &   18,225 & 6 & 1.0e$-$01 &  6,604 &  370,973 &   377,706 &  190.84 & 1.3e$-$04 & 11,639 &  49,969 &   30.57 \\         
 140 &   19,600 & 6 & 1.3e$-$01 &  6,616 &  302,593 &   309,288 &  171.16 & 1.4e$-$04 & 12,406 &  53,796 &   35.35 \\         
 145 &   21,025 & 6 & 1.5e$-$01 &  6,661 &  278,673 &   285,397 &  170.61 & 1.4e$-$04 & 12,814 &  56,423 &   40.04 \\         
 150 &   22,500 & 6 & 1.7e$-$01 &  7,048 &  269,053 &   276,151 &  178.26 & 1.5e$-$04 & 13,741 &  61,042 &   47.26 \\         
 155 &   24,025 & 6 & 1.9e$-$01 &  6,531 &  273,073 &   279,666 &  197.21 & 1.5e$-$04 & 14,651 &  65,802 &   54.88 \\         
 160 &   25,600 & 6 & 2.0e$-$01 &  6,204 &  261,354 &   267,616 &  203.83 & 1.6e$-$04 & 15,509 &  70,613 &   63.10 \\         
 165 &   27,225 & 6 & 2.1e$-$01 &  6,897 &  265,594 &   272,543 &  221.86 & 1.6e$-$04 & 15,310 &  69,910 &   64.23 \\         
 170 &   28,900 & 6 & 2.2e$-$01 &  5,214 &  286,714 &   292,023 &  257.46 & 1.7e$-$04 & 17,338 &  80,401 &   79.21 \\         
 175 &   30,625 & 6 & 2.6e$-$01 &  5,618 &  255,294 &   260,978 &  244.98 & 1.7e$-$04 & 18,015 &  83,984 &   87.44 \\         
 180 &   32,400 & 6 & 2.9e$-$01 &  5,260 &  237,754 &   243,072 &  242.00 & 1.8e$-$04 & 18,510 &  86,308 &   96.99 \\         
 185 &   34,225 & 6 & 3.1e$-$01 &  5,460 &  233,394 &   238,906 &  253.32 & 1.8e$-$04 & 19,457 &  91,745 &  114.16 \\         
 190 &   36,100 & 6 & 3.4e$-$01 &  5,741 &  217,814 &   223,590 &  259.29 & 1.9e$-$04 & 20,846 &  99,200 &  132.56 \\         
 195 &   38,025 & 6 & 3.5e$-$01 &  5,539 &  226,995 &   232,580 &  280.39 & 1.9e$-$04 & 20,135 &  97,440 &  138.37 \\         
 200 &   40,000 & 6 & 3.5e$-$01 &  4,969 &  241,335 &   246,370 &  311.54 & 2.0e$-$04 & 21,875 & 107,288 &  146.27 \\         
\hline
\end{tabular}}
\end{center}
\caption{Performances of NITSOL and Accelerated DF-SANE in the ``Driven Cavity
  Problem''.}
\label{tab3}
\end{table}

\begin{table}[ht!]
\begin{center}
\resizebox{0.85\textwidth}{!}{
\begin{tabular}{|cr|crrrr|crrr|}
\hline
\multirow{2}{*}{$n_p$} & \multicolumn{1}{c|}{\multirow{2}{*}{$n$}} &
\multicolumn{5}{c}{NITSOL (Newton-GMRES)} & \multicolumn{4}{|c|}{Accelerated DF-SANE} \\
\cline{3-11}
& &
\multicolumn{1}{c}{$\|F(x_*)\|_2$} &
\multicolumn{1}{c}{\#it$_1$} & \multicolumn{1}{c}{\#it$_{2}$} & \multicolumn{1}{c}{fcnt} &
\multicolumn{1}{c}{Time} & \multicolumn{1}{|c}{$\|F(x_*)\|_2$} & \multicolumn{1}{c}{\#it} &
\multicolumn{1}{c}{fcnt} & \multicolumn{1}{c|}{Time} \\
\hline
\hline
  100 &    10,000 & 7.9e$-$05 &     4 &      212 &      217 &    0.03 & 9.9e$-$05 &    174 &     349 &    0.15 \\         
  150 &    22,500 & 1.1e$-$04 &     4 &      321 &      326 &    0.14 & 1.5e$-$04 &    257 &     515 &    0.54 \\         
  200 &    40,000 & 1.6e$-$04 &     4 &      464 &      469 &    0.40 & 2.0e$-$04 &    335 &     671 &    1.42 \\         
  250 &    62,500 & 1.9e$-$04 &     4 &      655 &      660 &    0.98 & 2.5e$-$04 &    405 &     811 &    2.77 \\         
  300 &    90,000 & 2.4e$-$04 &     4 &      741 &      746 &    1.59 & 3.0e$-$04 &    473 &     947 &    4.95 \\         
  350 &   122,500 & 2.8e$-$04 &     4 &    1,054 &    1,059 &    3.15 & 3.4e$-$04 &    543 &   1,087 &    7.96 \\         
  400 &   160,000 & 3.2e$-$04 &     4 &    1,315 &    1,320 &    5.17 & 4.0e$-$04 &    613 &   1,227 &   12.47 \\         
  450 &   202,500 & 3.4e$-$04 &     4 &    1,518 &    1,523 &    7.54 & 4.3e$-$04 &    684 &   1,369 &   17.74 \\         
  500 &   250,000 & 4.1e$-$04 &     3 &    1,607 &    1,611 &   10.02 & 4.9e$-$04 &    754 &   1,509 &   25.21 \\         
  550 &   302,500 & 4.5e$-$04 &     3 &    2,447 &    2,451 &   19.94 & 5.5e$-$04 &    823 &   1,647 &   34.23 \\         
  600 &   360,000 & 5.0e$-$04 &     3 &    2,793 &    2,797 &   29.30 & 5.9e$-$04 &    893 &   1,787 &   43.80 \\         
  650 &   422,500 & 5.2e$-$04 &     4 &    3,502 &    3,507 &   45.42 & 6.3e$-$04 &    962 &   1,925 &   56.59 \\         
  700 &   490,000 & 5.6e$-$04 &     4 &    3,691 &    3,696 &   59.89 & 7.0e$-$04 &  1,030 &   2,061 &   71.45 \\         
  750 &   562,500 & 6.2e$-$04 &     4 &    3,702 &    3,707 &   72.01 & 7.4e$-$04 &  1,098 &   2,197 &   91.20 \\         
  800 &   640,000 & 6.4e$-$04 &     5 &    4,871 &    4,877 &  113.52 & 8.0e$-$04 &  1,164 &   2,329 &  111.31 \\         
  850 &   722,500 & 6.9e$-$04 &     5 &    4,826 &    4,832 &  132.01 & 8.5e$-$04 &  1,229 &   2,459 &  132.07 \\         
  900 &   810,000 & 7.2e$-$04 &     6 &    5,965 &    5,972 &  187.16 & 9.0e$-$04 &  1,292 &   2,585 &  157.64 \\         
  950 &   902,500 & 7.6e$-$04 &     6 &    5,695 &    5,702 &  200.12 & 9.5e$-$04 &  1,353 &   2,707 &  184.62 \\         
1,000 & 1,000,000 & 8.0e$-$04 &     8 &    7,386 &    7,395 &  292.89 & 1.0e$-$03 &  1,408 &   2,817 &  212.07 \\         
1,050 & 1,102,500 & 8.4e$-$04 &     7 &    6,901 &    6,909 &  304.10 & 1.0e$-$03 &  1,465 &   2,931 &  243.72 \\         
1,100 & 1,210,000 & 8.8e$-$04 &     9 &    8,399 &    8,409 &  405.89 & 1.1e$-$03 &  1,528 &   3,057 &  277.72 \\         
1,150 & 1,322,500 & 9.2e$-$04 &     8 &    7,977 &    7,986 &  423.56 & 1.1e$-$03 &  1,592 &   3,185 &  317.39 \\         
1,200 & 1,440,000 & 1.0e$-$03 &    10 &    9,980 &    9,991 &  578.89 & 1.2e$-$03 &  1,656 &   3,313 &  358.84 \\         
1,250 & 1,562,500 & 1.0e$-$03 &    10 &    9,678 &    9,689 &  608.95 & 1.2e$-$03 &  1,720 &   3,441 &  404.81 \\         
1,300 & 1,690,000 & 1.3e$-$03 &    11 &   10,960 &   10,972 &  747.45 & 1.3e$-$03 &  1,783 &   3,567 &  451.14 \\         
1,350 & 1,822,500 & 1.1e$-$03 &    11 &   10,935 &   10,947 &  803.39 & 1.3e$-$03 &  1,846 &   3,693 &  514.03 \\         
1,400 & 1,960,000 & 1.1e$-$03 &    13 &   12,932 &   12,946 & 1021.67 & 1.4e$-$03 &  1,909 &   3,819 &  618.58 \\         
1,450 & 2,102,500 & 1.2e$-$03 &    13 &   12,636 &   12,650 & 1074.93 & 1.4e$-$03 &  1,971 &   3,943 &  630.03 \\         
1,500 & 2,250,000 & 1.4e$-$03 &    14 &   13,920 &   13,935 & 1267.96 & 1.5e$-$03 &  2,032 &   4,065 &  701.65 \\         
\hline
\end{tabular}}
\end{center}
\caption{Performances of NITSOL and Accelerated DF-SANE in ``Generalized 2D
  Bratu Problem''.}
\label{tab4}
\end{table}

\begin{table}[ht!]
\begin{center}
\resizebox{0.85\textwidth}{!}{
\begin{tabular}{|cr|crrrr|crrr|}
\hline
\multirow{2}{*}{$n_p$} & \multicolumn{1}{c|}{\multirow{2}{*}{$n$}} &
\multicolumn{5}{c}{NITSOL (Newton-GMRES)} & \multicolumn{4}{|c|}{Accelerated DF-SANE} \\
\cline{3-11}
& &
\multicolumn{1}{c}{$\|F(x_*)\|_2$} &
\multicolumn{1}{c}{\#it$_1$} & \multicolumn{1}{c}{\#it$_{2}$} & \multicolumn{1}{c}{fcnt} &
\multicolumn{1}{c}{Time} & \multicolumn{1}{|c}{$\|F(x_*)\|_2$} & \multicolumn{1}{c}{\#it} &
\multicolumn{1}{c}{fcnt} & \multicolumn{1}{c|}{Time} \\
\hline
\hline
  100 &    10,000 & 7.9e$-$05 &    11 &    1,806 &    1,818 &    0.22 & 1.0e$-$04 &    684 &   1,376 &    0.52 \\         
  150 &    22,500 & 1.2e$-$04 &    10 &    2,666 &    2,677 &    0.92 & 1.5e$-$04 &    845 &   1,698 &    1.55 \\         
  200 &    40,000 & 1.7e$-$04 &    10 &    3,634 &    3,645 &    2.48 & 2.0e$-$04 &  1,255 &   2,521 &    4.68 \\         
  250 &    62,500 & 2.4e$-$04 &    12 &    5,817 &    5,830 &    6.83 & 2.5e$-$04 &  1,378 &   2,798 &    8.43 \\         
  300 &    90,000 & 2.6e$-$04 &    13 &    7,272 &    7,286 &   12.84 & 3.0e$-$04 &  1,595 &   3,212 &   14.81 \\         
  350 &   122,500 & 3.4e$-$04 &    13 &    7,890 &    7,904 &   19.29 & 3.5e$-$04 &  1,634 &   3,310 &   21.28 \\         
  400 &   160,000 & 3.9e$-$04 &    13 &    8,655 &    8,669 &   28.19 & 4.0e$-$04 &  1,616 &   3,349 &   29.48 \\         
  450 &   202,500 & 4.2e$-$04 &    15 &   10,911 &   10,927 &   45.61 & 4.5e$-$04 &  1,723 &   3,642 &   40.69 \\         
  500 &   250,000 & 4.7e$-$04 &    15 &   10,999 &   11,015 &   61.56 & 5.0e$-$04 &  1,858 &   3,916 &   55.15 \\         
  550 &   302,500 & 4.7e$-$04 &    18 &   14,065 &   14,084 &  103.30 & 5.5e$-$04 &  1,954 &   3,997 &   72.60 \\         
  600 &   360,000 & 5.3e$-$04 &    18 &   14,140 &   14,159 &  132.62 & 6.0e$-$04 &  2,043 &   4,145 &   90.35 \\         
  650 &   422,500 & 6.3e$-$04 &    17 &   13,152 &   13,170 &  151.83 & 6.5e$-$04 &  2,171 &   4,508 &  116.74 \\         
  700 &   490,000 & 6.4e$-$04 &    20 &   16,127 &   16,148 &  229.96 & 7.0e$-$04 &  2,271 &   4,874 &  144.36 \\         
  750 &   562,500 & 7.3e$-$04 &    18 &   14,164 &   14,183 &  243.82 & 7.5e$-$04 &  2,377 &   5,263 &  180.96 \\         
  800 &   640,000 & 7.2e$-$04 &    20 &   16,127 &   16,148 &  344.09 & 8.0e$-$04 &  2,502 &   5,630 &  220.13 \\         
  850 &   722,500 & 7.6e$-$04 &    20 &   16,147 &   16,168 &  396.81 & 8.5e$-$04 &  2,629 &   6,037 &  263.11 \\         
  900 &   810,000 & 8.8e$-$04 &    20 &   16,129 &   16,150 &  431.85 & 9.0e$-$04 &  2,750 &   6,398 &  309.54 \\         
  950 &   902,500 & 9.3e$-$04 &    21 &   17,094 &   17,116 &  521.95 & 9.5e$-$04 &  2,909 &   6,856 &  367.30 \\         
1,000 & 1,000,000 & 9.7e$-$04 &    20 &   16,144 &   16,165 &  549.77 & 1.0e$-$03 &  2,986 &   7,227 &  419.13 \\         
1,050 & 1,102,500 & 9.9e$-$04 &    22 &   18,098 &   18,121 &  679.61 & 1.0e$-$03 &  3,048 &   7,581 &  489.69 \\         
1,100 & 1,210,000 & 1.0e$-$03 &    22 &   18,111 &   18,134 &  753.17 & 1.1e$-$03 &  3,125 &   8,017 &  548.80 \\         
1,150 & 1,322,500 & 1.1e$-$03 &    22 &   18,026 &   18,049 &  892.35 & 1.1e$-$03 &  3,264 &   8,448 &  631.76 \\         
1,200 & 1,440,000 & 1.2e$-$03 &    23 &   19,006 &   19,030 & 1056.51 & 1.2e$-$03 &  3,434 &   9,051 &  748.79 \\         
1,250 & 1,562,500 & 1.1e$-$03 &    24 &   19,983 &   20,008 & 1183.49 & 1.2e$-$03 &  3,371 &   9,061 &  815.04 \\         
1,300 & 1,690,000 & 1.3e$-$03 &    24 &   20,010 &   20,035 & 1207.67 & 1.3e$-$03 &  3,377 &   9,249 &  891.39 \\         
1,350 & 1,822,500 & 1.3e$-$03 &    26 &   21,969 &   21,996 & 1440.00 & 1.3e$-$03 &  3,302 &   8,892 &  875.23 \\         
1,400 & 1,960,000 & 1.3e$-$03 &    26 &   21,947 &   21,974 & 1538.56 & 1.4e$-$03 &  3,234 &   8,477 &  893.25 \\         
1,450 & 2,102,500 & 1.4e$-$03 &    27 &   22,938 &   22,966 & 2118.16 & 1.4e$-$03 &  3,284 &   8,566 &  978.04 \\         
1,500 & 2,250,000 & 1.4e$-$03 &    27 &   22,934 &   22,962 & 2108.73 & 1.5e$-$03 &  3,314 &   8,569 & 1047.10 \\         
\hline
\end{tabular}}
\end{center}
\caption{Performances of NITSOL and Accelerated DF-SANE in ``Flow in a Porous
  Media Problem''.}
\label{tab5}
\end{table}

\subsection{Comparison with Anderson Mixing}

In a third set of experiments, we compared Accelerated DF-SANE with
Anderson Mixing, given by (\ref{anmix}--\ref{anmix3}) and described
in~\cite[\S 2.5]{fangsaad}. According to~\cite{eyert}
and~\cite{fangsaad}, Andrerson Mixing is equivalent to a
generalization of Broyden's second method. Once again, we do not
consider preconditioners in the experiments, since they could be used
in all problems regardless of the method considered. Anderson Mixing
was implemented in Fortran~90. In the implementation, the minimum-norm
solution to the least-squares problems is computed with th same
updated QR decomposition that is used in the Accelerated DF-SANE
method, as described at the end of Section~\ref{secQR}.

\begin{figure}[ht!]
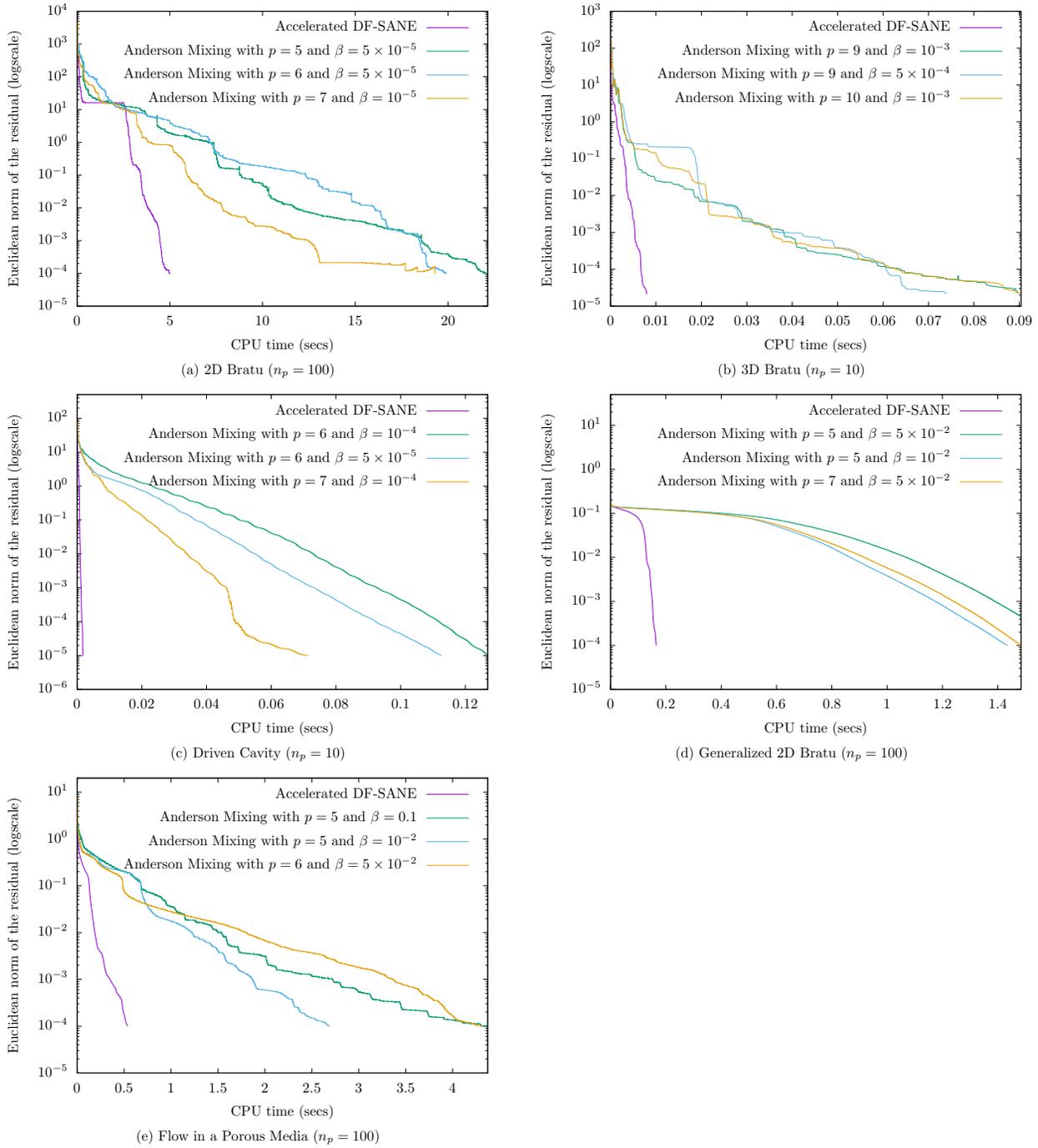

\begin{center}
\resizebox{\textwidth}{!}{\begin{tabular}{cc}
\input{bmdiisfig41.tex} & \input{bmdiisfig42.tex} \\
(a) 2D Bratu ($n_p=100$) & (b) 3D Bratu ($n_p=10$) \\
\input{bmdiisfig43.tex} & \input{bmdiisfig44.tex} \\
(c) Driven Cavity ($n_p=10$) & (d) Generalized 2D Bratu ($n_p=100$) \\
\input{bmdiisfig45.tex} &  \\
(e) Flow in a Porous Media ($n_p=100$) & \\
\end{tabular}}
\end{center}
\caption{Comparison between Accelerated DF-SANE and Anderson Mixing.}
\label{fig4}
\end{figure}

The experiments in~\cite{fangsaad} make it clear that the performance
of Anderson Mixing depends strongly on the choice of the $\beta$
parameter; and that this parameter varies from problem to problem and
even when the same problem with different dimensions is solved; see,
for example, \cite[\S 5.1]{fangsaad} where a variant of Bratu problem
is solved. As a consequence, it would be impracticable to address each
of the five problems in Sections~\ref{sec61} and~\ref{sec62} by
varying the values of $n_p$ as we did in
Tables~\ref{tab1}--\ref{tab5}. In fact we tried to do that but, for
example, if in the 2D Bratu problem with $n_p=125$ we use the best
value of~$p$ and~$\beta$ that we found for the case $n_p=100$, the
method diverges. That is an issue when one wants to do experiments
like the ones in the previous sections, but it is not a serious
problem if one has a practical problem to solve. Therefore, we
considered in this comparison only the first problem of each of
Tables~\ref{tab1}--\ref{tab5}. For each of the problems, we made an
exhaustive search for the best values of $p$ and $\beta$ which, based
on the values of $\beta$ reported in~\cite{fangsaad}, comprised
testing all combinations of $p \in \{5,6,\dots,10\}$ and $\beta \in \{
\beta_a 10^{-\beta_b} \; | \; \beta_a \in \{0.5, 1\} \mbox{ and }
\beta_b=1,\dots,8 \}$, totalizing 96 combinations. Figure~\ref{fig4}
show, for each of the five problems considered, the evolution of the
residual norm over time obtained with the three best (fastest) combinations of
parameters. (Combinations are different for each problem.) The figures
also include the results obtained with Accelerated DF-SANE. The result
shown for Accelerated DF-SANE is the one already reported in the
previous sections, without any extra tuning of its parameters. The
figure shows that Accelerated DF-SANE is approximately five to ten
times faster than Anderson Mixing in the considered problems. The
comparison is essentially the same if the number of function
evaluations is compared.

\section{Conclusions} \label{conclusions}

The Sequential Secant Method, which is the most obvious
multidimensional generalization of the secant method for solving
nonlinear equations, seems to have been introduced by Wolfe in 1959
\cite{wolfe}; see also \cite{barnes}. This method has been analyzed in
classical books and papers \cite{or,schwetlick,jankowska} where, under
suitable conditions, local convergence with $R$-order equal to the
unique solution of $t^{n+1} - t^n - 1 = 0$ was proved. Given the
consecutive iterates $x^{k-n}, \dots, x^{k-1}, x^k$, the Sequential
Secant Method computes
\begin{equation} \label{secmas}
x^{k+1} = x^k - (s^{k-n}, \dots, s^{k-1})(y^{k-1}, \dots,
y^{k-1})^{-1}F(x^k).
\end{equation}
This iteration is well defined if the matrix $(y^{k-n}, \dots,
y^{k-1})$ is nonsingular.  Moreover the good local convergence
properties need uniformly linear independence of the increments
$s^{k-n}, \dots, s^{k-1}$. The method has been updated in several ways
in order to fix these drawbacks while maintaining its convergence
properties. The natural limited memory version of (\ref{secmas}) is
defined by
\begin{equation} \label{seclim}
x^{k+1} = x^k - (s^{k-p}, \dots, s^{k-1})(y^{k-p}, \dots,
y^{k-1})^{\dagger}F(x^k),
\end{equation}   
where $1 \leq p \ll n$. This formula is inconvenient for solving
nonlinear systems because $s^k$ necessarily belongs to the subspace
generated by $s^{k-p}, \dots, s^{k-1}$ which implies that $x^{k+j}$ is
in the affine subspace determined by $x^{k-p}, \dots, x^{k-1}, x^k$,
for all $j$, and, so, convergence to a solution cannot occur unless
the solution belongs to the same affine subspace. This is the reason
why, in the present paper, we do not use the method defined by
(\ref{seclim}). Instead, we compute $x^{k+1}_{\trial}$ using the
Sequential Residual approach, we define $s^k = x^{k+1}_{\trial} -
x^k$, $y^k=F(x^{k+1}_{\trial}) - F(x^k)$,
\begin{equation} \label{xacel}
x^{k+1}_{\acc} =  x^{k+1}_{\trial} - (s^{k-p+1}, \dots, s^{k-1}, s^k)(y^{k-p+1},
\dots, y^{k-1}, y^k)^{\dagger}F(x^{k+1}_{\trial}),
\end{equation}
and we choose $x^{k+1}$ as the best of the trials $x^{k+1}_{\trial}$
and $x^{k+1}_{\acc}$. In this way, we preserve the good properties of
the Sequential Secant Method associated with the good global behavior
of Sequential Residual approaches. The freedom on the choice of the
residual step favors the employment of preconditioners when they are
available.

The most popular methods for solving nonlinear systems of equations in
which the application of Newton's method is impossible or extremely
expensive are based on the Inexact Newton approach with
Krylov-subspace methods (as GMRES) for solving approximately the
Newtonian linear systems at each iteration. These methods have a long
tradition and, very likely, they deserve to be the preferred ones by
practitioners in the numerical PDE community. Nevertheless, in this
paper we showed that, for some very large interesting problems, an
approach based on sequential-secant-like accelerations of a residual
method is more efficient than a standard implementation of
Newton-GMRES with its default algorithmic parameters.  This indicates
that those problems possess characteristics that favor the application
of the secant paradigm over the inexact Newton one. Of course, the
opposite situation probably occurs in many cases. This means,
therefore, that efficiency for solving practical problems will be
increased if practitioners have easy access to both types of
methods. The codes used in this paper, that make it possible the
reproducibility of all the experiments, as well as the instructions
for using our algorithm may be found in
\url{http://www.ime.usp.br/~egbirgin/sources/accelerated-df-sane/}.

Future work will include the employment of the new methods to the
acceleration of KKT solvers that are necessary in the Augmented
Lagrangian approach for constrained optimization~\cite{bmbook}.\\

\noindent
\textbf{Acknowledgments.} The authors would like to thank the referees
for their careful reading and the many suggestions they contributed to
improve the quality of this work.

\end{document}